\documentclass[twoside,leqno,10pt]{article}

\usepackage{graphics,ifthen}
\usepackage{latexsym}
\usepackage{amsmath}
\usepackage{amssymb}
\usepackage{mathrsfs}
\begin{document}
\input{latexP.sty}
\input {epsf.sty}
\def\ind{\stackrel{\mathrm{ind}}{\sim}}
\def\iid{\stackrel{\mathrm{iid}}{\sim}}
\def\Prodi{\mathop{{\lower9pt\hbox{\epsfxsize=15pt\epsfbox{pi.ps}}}}}
\def\prodi{\mathop{{\lower3pt\hbox{\epsfxsize=7pt\epsfbox{pi.ps}}}}}
\def\Definition{\stepcounter{definitionN}
    \Demo{Definition\hskip\smallindent\thedefinitionN}}
\def\EndDefinition{\EndDemo}
\def\Example#1{\Demo{Example [{\rm #1}]}}
\def\EndExample{\qed\EndDemo}
\def\Category#1{\centerline{\Heading #1}\rm}
\def\e{\text{\hskip1.5pt e}}
\newcommand{\eps}{\epsilon}
\newcommand{\proof}{\noindent {\bf Proof:\ }}
\newcommand{\remarks}{\noindent {\bf Remarks:\ }}
\newcommand{\note}{\noindent {\bf Note:\ }}
\newcommand{\examp}{\noindent {\bf Example:\ }}
\newcommand{\Lower}[2]{\smash{\lower #1 \hbox{#2}}}
\newcommand{\ben}{\begin{enumerate}}
\newcommand{\een}{\end{enumerate}}
\newcommand{\bi}{\begin{itemize}}
\newcommand{\ei}{\end{itemize}}
\newcommand{\hp}{\hspace{.2in}}

\newcommand{\Bcr}{\mathscr{B}}
\newcommand{\Ucr}{\mathscr{U}}
\newcommand{\Gcr}{\mathscr{G}}
\newcommand{\Dcr}{\mathscr{D}}
\newcommand{\CS}{\mathscr{C}}
\newcommand{\Fcr}{\mathscr{F}}
\newcommand{\Icr}{\mathscr{I}}
\newcommand{\Lcr}{\mathscr{L}}
\newcommand{\Mcr}{\mathscr{M}}
\newcommand{\Ncr}{\mathscr{N}}
\newcommand{\Pcr}{\mathscr{P}}
\newcommand{\Qcr}{\mathscr{Q}}
\newcommand{\Scr}{\mathscr{S}}
\newcommand{\Tcr}{\mathscr{T}}
\newcommand{\Xcr}{\mathscr{X}}
\newcommand{\E}{\mathbb{E}}
\newcommand{\F}{\mathbb{F}}
\newcommand{\I}{\mathbb{I}}
\newcommand{\Q}{\mathbb{Q}}
\newcommand{\X}{\mathbb{X}}
\newcommand{\Pe}{\mathbb{P}}
\newcommand{\M}{\mathbb{M}}
\newcommand{\R}{\mathbb{R}}
\newcommand{\Wbb}{\mathbb{W}}

\newcommand{\Xn}{X_1,\ldots,X_n}
\newcommand{\lan}{\langle}
\newcommand{\ra}{\rangle}
\newcommand{\edr}{\mathrm{e}}
\newcommand{\ddr}{\mathrm{d}}
\newcommand{\Var}{\operatornamewithlimits{Var}}
\newcommand{\ld}{\ldots}
\newcommand{\cd}{\cdots}
\def\simind{\stackrel{\mbox{\scriptsize{ind}}}{\sim}}
\def\simiid{\stackrel{\mbox{\scriptsize{iid}}}{\sim}}

\newcommand{\Ga}{\Gamma}
\newcommand{\si}{\sigma}
\newcommand{\la}{\lambda}
\newcommand{\ph}{\varphi}
\newcommand{\ep}{\varepsilon}
\newcommand{\ga}{\gamma}
\newcommand{\Om}{\Omega}
\newcommand{\om}{\omega}

\newtheorem{thm}{Theorem}[section]
\newtheorem{defin}{Definition}[section]
\newtheorem{prop}{Proposition}[section]
\newtheorem{lem}{Lemma}[section]
\newtheorem{cor}{Corollary}[section]
\newcommand{\rb}[1]{\raisebox{1.5ex}[0pt]{#1}}
\newcommand{\mc}{\multicolumn}
\def\Beta{\text{Beta}}
\def\Dir{\text{Dirichlet}}
\def\DP{\text{DP}}
\def\P{{\bf p}}
\def\fhat{\widehat{f}}
\def\GA{\text{gamma}}
\def\ind{\stackrel{\mathrm{ind}}{\sim}}
\def\iid{\stackrel{\mathrm{iid}}{\sim}}
\def\J{{\bf J}}
\def\K{{\bf K}}
\def\min{\text{min}}
\def\N{\text{N}}
\def\p{{\bf p}}
\def\Sj{{\bf S}}
\def\sj{{\bf s}}
\def\U{{\bf U}}
\def\u{{\bf u}}
\def\w{{\bf w}}
\def\W{{\bf W}}
\def\X{{\bf X}}
\def\x{{\bf x}}
\def\y{{\bf y}}
\def\Y{{\bf Y}}
\newcommand{\reals}{{\rm I\!R}}
\newcommand{\PR}{{\rm I\!P}}
\def\Z{{\bf Z}}
\def\yy{{\mathcal Y}}
\def\rr{{\mathcal R}}
\def\mm{{\mathcal M}}
\def\BP{\text{beta}}
\def\ts{\tilde{t}}
\def\Ns{\tilde{N}}
\def\Ps{\tilde{P}}
\def\gs{\tilde{g}}
\def\fs{\tilde{f}}
\def\ys{\tilde{Y}}
\def\ps{\tilde{P}}
\def
\Report{\centerline{\small{\rm NRM Bayesian Models}}}
\def\Author{\centerline{\small{\rm L.F. JAMES, A. LIJOI AND I.
PR\"UNSTER}}} \pagestyle{myheadings} \markboth{\Author}{\Report}
\thispagestyle{empty} \bct\Heading BAYESIAN INFERENCE VIA CLASSES
OF\\ NORMALIZED RANDOM MEASURES \lbk\lbk\smc {\sc Lancelot F. James,
Antonio Lijoi and Igor Pr\"unster}\footnote{ \eightit AMS 2000
subject classifications.
              \rm Primary 62G05; secondary 62F15.\\
\indent\eightit Keywords and phrases. \rm
         Bayesian Nonparametrics,
         Chinese restaurant process,
         Generalized gamma convolutions,
         Gibbs partitions,
         Poisson random measure
         }
\lbk\lbk \BigSlant Hong Kong University of Science and Technology
and University of Pavia\rm \lbk 
\ect \Quote One of the main research areas in Bayesian
Nonparametrics is the proposal and study of priors which
generalize the Dirichlet process. Here we exploit theoretical
properties of Poisson random measures in order to provide a
comprehensive Bayesian analysis of random probabilities which are
obtained by an appropriate normalization. Specifically we achieve
explicit and tractable forms of the posterior and the marginal
distributions, including an explicit and easily used description
of generalizations of the important Blackwell-MacQueen P\'olya urn
distribution. Such simplifications are achieved by the use of a
latent variable which admits quite interesting interpretations
which allow to gain a better understanding of the behaviour of
these random probability measures. It is noteworthy that these
models are generalizations of models considered by Kingman~(1975)
in a non-Bayesian context. Such models are known to play a
significant role in a variety of applications including genetics,
physics, and work involving random mappings and assemblies. Hence
our analysis is of utility in those contexts as well. We also show
how our results may be applied to Bayesian mixture models and
describe computational schemes which are generalizations of known
efficient methods for the case of the Dirichlet process. We
illustrate new examples of processes which can play the role of
priors for Bayesian nonparametric inference and finally point out
some interesting connections with the theory of generalized gamma
convolutions initiated by Thorin and further developed by
Bondesson.\EndQuote

\rm \section{Introduction}  A key problem in Bayesian nonparametric
inference is the definition of a prior distribution on the space of
all probability measures. Starting from the papers by Ferguson
(1973) and Freedman (1963), in which the celebrated Dirichlet
process has been introduced, various approaches for constructing
random probability measures, whose distribution acts as a
nonparametric prior, have been undertaken. They all aim at providing
generalizations of the Dirichlet process. Among them we mention the
neutral-to-the-right random probability measures due to Doksum
(1974), which are obtained via an exponential transformation of an
increasing process with independent increments, and the
P\'olya--tree--priors thoroughly studied by Mauldin, Sudderth and
Williams (1992) and Lavine (1992), which arise by considering
suitable urn schemes on trees of nested partitions. Moreover, it is
well-known that the Dirichlet process can be defined by normalizing
the increments of a gamma process. The idea of constructing random
probability measures by means of a normalization procedure has been
exploited and developed in a variety of contexts not closely related
to Bayesian inference. Indeed it has found many interesting
applications: Kingman (1975) and Janson (2001) for storage problems
and applications to computer science; Ewens and Tavar\'e (1995) and
Grote and Speed (2002) for population genetics; Engen (1978) and
McCloskey (1965) for ecology; Derrida (1981) and Ruelle (1987) for
statistical physics; Donnelly and Grimmett (1993) and Pitman (2002)
for combinatorics and number theory; Pitman (1997) and Pitman and
Yor (1997) for excursion theory.

Kingman's~(1975) paper suggests that one can construct random
probability measures as follows. First take the ranked points of a
homogeneous Poisson process on $(0,\infty)$, say
$\Delta_{1}>\Delta_{2},\ldots$, such that their sum
$\sum_{i=1}^{\infty}\Delta_{i}$ is finite and positive almost
surely. Use these points to construct a sequence of probabilities,
$(Q_{i}=\Delta_{i}/\sum_{l=1}^{\infty}\Delta_{l})$. Independent of
this sequence choose $(Z_{i})$ to be an iid sequence of random
elements of a Polish space with common distribution, say $H$. A
random probability measure is then formed by
$\sum_{i=1}^{\infty}Q_{i}\delta_{Z_{i}}$. Two special cases include
the Dirichlet process [Ferguson~(1973)] and a class of random
probability measures based on the stable law discussed in for
instance Pitman~(1996). It is interesting to note that much of the
analysis related to such quantities focuses exclusively on the
behaviour of the sequence $(Q_{i})$. Such studies were carried out
for instance by Pitman, Perman and Yor~(1992) and Pitman~(2003).
Note that because of the independence of $(Q_{i})$ and $(Z_{i})$,
the distinctive features of different classes of random probability
measures are in fact deduced exclusively from the analysis of
$Q_{i}$. However with the exception of the Dirichlet process and
those models based on a stable law, the analysis related to the bulk
of these processes has yet to yield tractable results suitable for
practical implementation in Bayesian nonparametric problems. This is
due in part to the fact that the focus in Kingman~(1975) and the
majority of the subsequent analysis involves considerations other
than Bayesian applications. Thus the issues of general tractability,
in terms of their possible usage in a Bayesian context, raised by
Adrian Smith and others in Kingman (1975), remains open.

In this paper we provide an analysis, for a larger class of
models, with a view toward practical implementation and a better
theoretical understanding of such models in a Bayesian context.
Using the Dirichlet process as a benchmark, such an analysis
requires a suitable description of the posterior distribution,
analogous to that given by Ferguson~(1973). Additionally, noting
the success of the Dirichlet process in complex mixture models, we
find tractable analogues of the Blackwell-MacQueen P\'olya Urn and
accompanying Chinese restaurant distribution which is related to
the Ewens sampling formula derived by Ewens~(1972) and Antoniak
(1974). This then paves the way for the description of
implementable MCMC and SIS computational procedures to approximate
efficiently posterior quantities in applications such as
hierarchical mixture modeling. Also, the exchangeable marginal
distribution is equivalent to the notion of the moment measures of
a random probability measure and hence is basic to the
understanding of its theoretical properties and to the calculation
of higher order moments. Of course, as a natural byproduct, our
analysis has implications in the related non-Bayesian contexts
described above. Our methodology follows the paradigm laid out in
James~(2005a, 2002). James (2005a) points out that, in analogy to
the use of classical Bayes rule, one often has to introduce
additional refinements to obtain the most tractable forms. Here we
show that such refinements are achieved by the use of an important
exploitation of a {\it latent} structure derived from the gamma
identity.

We shall actually analyze a richer and more complex generalization
of the model suggested by Kingman~(1975), which have recently
arisen in a Bayesian context. Regazzini, Lijoi and Pr\"unster
(2003) consider random probability measures obtained by a
normalization of suitably time-changed increasing processes with
independent and not necessarily stationary increments. Their
interest was primarily the formidable problem of the determining
the distribution of mean functionals. That is a generalization of
the important body of work initiated by Cifarelli and
Regazzini~(1990).  James~(2002), using an approach closely
connected to Perman, Pitman and Yor~(1992), considers a more
general $h$-biased variation of the random probability measures
which allows for an extension to arbitrary Polish spaces. The
constructions coincide on Euclidean spaces when $h(s)=s$, for
$s\in(0,\infty)$. In any case, such models can also be represented
in terms of corresponding $(Q_{i})$ and $(Z_{i})$. However, now
these two sequences are not necessarily independent, meaning one
cannot just analyze the $(Q_{i})$. Moreover the $(Q_{i})$ are now
generated by an $h$-biased structure.

We also provide a more specific analysis of three classes of random
probability measures which are somehow connected to the Dirichlet
process, and hence inherit many of its desirable features, but are
otherwise much more general. Specifically we define and examine a
class of dependent Dirichlet processes, a class of random
probability connected to the beta process. Lastly, using the
$h$-biased framework,  we construct random probability measures
based on the theory of Generalized Gamma Convolutions (GGC),
initiated by Thorin~(1977, 1978) and provide a detailed analysis.
Some interesting features of this construction, is that we are able
to embed a very large class of random probability measures within a
tractable framework of models which are extensions of Dirichlet
process mean functionals. More precisely, these are models derived
from a gamma process with a shape measure which is sigma finite and
therefore has possibly infinite total mass. Our construction
includes, for instance the stable law random probability measure and
hence the entire class of the two-parameter Poisson-Dirichlet
process. Moreover our results and discussions are strongly
suggestive of a more synergistic interplay between the theory of GGC
and the study of the laws of extensions of mean functionals of the
Dirichlet process.

We next describe the formal construction and discuss the key
features of our analysis.

\subsection{Preliminaries}
Let $N$ denote a Poisson random measure on an arbitrary Polish
Space $\Scr \times\Xcr$, with mean intensity
$$
\E\[N(ds,dx)\]=\nu(ds,dx).$$ Denote the law of $N$ as
$\Pe(\cdot|\nu)$ with $\E[\,\cdot\,|\nu]$ denoting expectation
with respect to $\Pe$. The distribution of $N$ is characterized by
the Laplace functional
    \Eq
    \label{loglaplace}
    \E\left[\edr^{-N(g)}\right]=\exp\(-\int_{\Xcr\times\Scr}(1-\edr^{-g(s,x)})\,\nu(\ddr
    s,\, \ddr x)\)
    \EndEq
where $N(g)=\int_{\Scr\times\Xcr} g(s,x)N(ds,dx)$ and $g$ is any
positive function. Throughout, as in Daley and Vere-Jones~(1986),
we have that $N$ and all related functionals take their values in
the space of boundedly finite measures on $\Scr\times \Xcr$, say
$\M$. A measure, say $\Lambda$, is boundedly finite if for each
bounded set $A$, $\Lambda(A)<\infty.$ According to the
decomposition of the intensity measure, we distinguish the two
following cases: \Enumerate
\item[(i)] if $\nu(\ddr s,\, \ddr x)=\rho(\ddr s)\,H(\ddr
x)$ we say that $N$ and its related functionals are
\textit{homogeneous}\\[-15pt]
\item[(ii)] if $\nu(\ddr s,\, \ddr x)=\rho(\ddr s|x)\,\eta(\ddr
x)=F_{0}(\ddr x|s)\rho(\ddr s)$, $N$ and its related functionals are
\textit{non-homogeneous}. \EndEnumerate Where $H$ and $F_{0}$ are
probabilities on $\Xcr$, $\eta$ is a finite measure on $\Xcr$. The
quantity, $\rho(ds|x)$ is an inhomogeneous L\'evy measure for each
fixed $x$. The quantity,
$$
\rho(ds)=\int_{\Xcr}\rho(ds|x)\eta(dx)
$$
is a then a homogeneous L\'evy measure on $\Scr$. Let
$(J_{i},Z_{i})$ denote the sequence of points of $N$ on
$\Scr\times\Xcr$. Using the decomposition $F_{0}(dx|s)\rho(ds)$,
it follows that now for each $i$, the conditional distribution of
$Z_{i}$ given $J_{i}$ is
$$
\PR(Z_{i}\in dx|J_{1},J_{2},\ldots):=\PR(Z_{i}\in
dx|J_{i})=F_{0}(dx|J_{i}).$$ Moreover the marginal distribution of
the sequence $(J_{i})$ are the random points in $\Scr$ of a
homogeneous Poisson random measure with mean intensity $\rho(ds)$.
See for instance Resnick~(1987, section 3.3.2).

Let $h$ denote a strictly positive function on $\Scr$. Now define
a random measure on $\Xcr$, which is representable in distribution
as,
$$
\mu(dx)=\int_{\Scr}h(s)N(ds,dx)
$$
such that its total mass $T:=\mu(\Xcr)=\sum_{i=1}^{\infty}h(J_{i})$,
is strictly positive and almost surely finite. This happens if
$\nu(\Scr\times\Xcr)=+\infty$ and the Laplace transform
    $$
\phi(\lambda)=\E\left[\edr^{-\lambda\,T}\right]={\mbox
e}^{-\psi(\lambda)}
$$
where
$$
 \psi(\lambda):=\int_{\Xcr\times\Scr}(1-\edr^{-\lambda\,
    h(s)})\,\nu(\ddr s,\, \ddr x) 
$$
is finite for any positive $\lambda$. We now can define a class of
random probability measures on $\Xcr$ representable in distribution
as, \Eq
P(dx)=\frac{\int_{\Scr}h(s)N(ds,dx)}{T}=\frac{\mu(dx)}{T}=\sum_{i=1}^{\infty}Q_{i}\delta_{Z_{i}}
=\sum_{i=1}^{\infty}\frac{h(J_{i})}{T}\delta_{Z_{i}} \label{NRM}
\EndEq We will call the class of $P$ defined as~\mref{NRM},
normalized random measures or NRMs. In Regazzini~et~al. (2003) an
analogous construction, with $\Scr=(0,+\infty)$ and $\Xcr=\R$, has
been developed by normalizing the increments of an increasing
additive process.

\Remark Note that the generality of $\Scr$ allows for quite
complex spaces. For example, one could take $\Scr$ to be the space
of probability measures. Or $\Scr$ could denote the space
containing the sample paths of Brownian Motion or more general
stochastic processes. See Perman, Pitman and Yor~(1992) for an
example involving excursion spaces of Markov processes. \EndRemark

\Remark Pitman~(2003) provides an important extension of the class
of homogeneous NRM, defined by  $h(s)=s\in(0,\infty)$, which can
be defined as follows. For a homogeneous $\rho$, denote the law of
the sequence $(J_{i}/T)$ as $PK(\rho)$ Furthermore, define the law
of $(J_{i}/T)|T=t$, as $PK(\rho|t)$,  where $T$ has density
$f_{T}(t)$. Then this class generates the class of laws on
$(J_{i}/T)$, given by
$PK(\rho;\gamma_{T})=\int_{0}^{\infty}PK(\rho|t)\gamma_{T}(dt)$,
where $\gamma_{T}$ is some density of $T$. An interesting special
case is the choice of $\gamma_{T}(t)\propto t^{-q}f_{T}(t)$,
which, with the exception of the Dirichlet process,  yields the
two-parameter Poisson Dirichlet process when $f_{T}$ is the
density of the stable law of index $0<\alpha<1$. The idea of
conditioning on $T$ is discussed briefly in Kingman~(1975). In
order to capture models like this and more general processes we
can instead work with a weighted Poisson law of the type
$g(N)\Pe(dN|\nu)$ where we assume that $\int_{\M}
g(N)\Pe(dN|\nu)=1$. Formal details are given in section 8.3.
\EndRemark

\section{Posterior Analysis} Similar to Ferguson~(1973) we shall
consider the classical setup. Let $(X_n)_{n\ge 1}$ be a sequence
of $\Xcr$--valued exchangeable random elements, i.e. such that for
any $n\ge 1$ the $X_{1},\ldots, X_{n}$ are, conditional on $P$,
iid with common distribution $P$. Suppose, also, that $P$ is an
NRM. This yields a description of the joint distribution of
$(\X,P)$. We are interested in obtaining its description in terms
of a posterior distribution of $P|\X$ and the (exchangeable)
marginal distribution of $\X$, say $\Mcr$. Since the law of $P$ is
dominated by the law of $N$ it follows that the posterior
distribution of $P|\X$ is determined by the posterior distribution
of $N|\X$. Moreover $\Mcr$ can be expressed as,$$
\Mcr(dx_{1},\ldots,dx_{n})=\int_{\M}\[\prod_{i=1}^{n}P(dx_{i})\]
\Pe(dN|\nu)
=\int_{\M}T^{-n}\[\prod_{i=1}^{n}\int_{\Scr}h(s_{i})N(ds_{i},dx_{i})\]
\Pe(dN|\nu)
$$
and is the general analogue of the Blackwell-MacQueen P\'olya Urn
distribution. We recall from James~(2005a) that there is a one to
one correspondence between $\X$ and $(\Y, \p)$, where using
notation similar to Lo~(1984), $\Y=(Y_1, \ldots, Y_{n(\p)})$
denotes the distinct values of $\X$ and $\p=\{C_{1},\ldots,
C_{n(\p)}\}$ stands for a partition of the integers $\{1, \ldots ,
n\}$ of size $n(\p)\leq n$ recording which observations are equal.
The number of elements in the $j$-th cell,
$C_{j}:=\{i:X_{i}=Y_{j}\}$, of the partition is indicated by
$e_{j}$, for $j=1, \ldots , n(\p)$, so that
$\sum_{j=1}^{n(\p)}e_{j}=n$. When it is necessary to emphasize a
further dependence on $n$, we will also use the notation
$e_{j,n}:=e_{j}$. It follows that the marginal distribution of
$\X$ can be expressed in terms of a conditional distribution of
$\X|\p$, which is the same as a conditional distribution of the
unique values $\Y|\p$, and the marginal distribution of $\p$. The
marginal distribution of $\p$, denoted as $\pi(\p)$ or
$p(e_{1},\ldots, e_{n(\p)})$, is an {\it exchangeable partition
probability function} (EPPF). \Remark A detailed general
discussion of the EPPF concept may be found in Pitman~(2002). A
discussion of its role in a Bayesian context pertinent to our
homogeneous models, which are a special case of species sampling
models, may be deduced from Pitman~(1996) and Ishwaran and
James~(2003a). Its role for inhomogeneous models is a bit
different. See James~(2005a) for the role of the EPPF, and that of
$\Mcr$, in this more general Bayesian context. The most well-known
EPPF is the Chinese restaurant distribution associated with the
Dirichlet process. \EndRemark

\subsection{Posterior Distributions}
Now, similar to Perman, Pitman and Yor~(1992, section 4) let,
$N_{n(\p)}:=N-\sum_{j=1}^{n(\p)}\delta_{J_{j,n},Y_{j}}$ denote a
random measure after the first $n(\p)$ pairs of unique points
$(J_{j,n},Y_{j})$ in $\Scr\times\Xcr$ are picked from $N$ by
$h$-biased sampling. Define $N_{0}:=N$ and note that the law of
$N_{n(\p)}$, depends on the random variable $n(\p)$. For each $n$ it
follows that $N=N_{n(\p)}+\sum_{j=1}^{n(\p)}\delta_{J_{j,n},Y_{j}}$.
From this one may define
$$
\mu_{n(\p)}=\mu-\sum_{j=1}^{n(\p)}h(J_{j,n})\delta_{Y_{j}}\quad
{\mbox { and }}\quad
T_{n(\p)}=\mu_{n(\p)}(\Xcr)=T-\sum_{j=1}^{n(\p)}h(J_{j,n})
$$
Now, crucial to our exposition, for each $n$ define the random
variable $U_{n}=\Gamma_{n}/T$, where $\Gamma_{n}$ denotes a gamma
random variable with shape $n$ and scale $1$ which is independent of
$N$ and hence $T$. Throughout our exposition $\Gamma_{n}$ will
always be taken to be independent of every random variable except
for $U_{n}$. It will be shown in the appendix, that the appearance
of $U_{n}$ arises from an application of the Gamma identity, \Eq
T^{-n}=\frac{1}{\Gamma(n)}\int_{0}^{\infty}u^{n-1}{\mbox e}^{-uT}du
\label{gammaid}\EndEq It turns out that importantly, $U_{n}$ is also
intimately connected to the {\it real inversion formula for the
Laplace transform}. This is discussed in Feller (1971, VII.6) and
plays a prominent role in the analysis of infinite-divisibility by
Thorin~(1977,1978) and Bondesson~(1979, 1992). We will discuss some
features of this in the forthcoming sections.

Define for each integer $l$, {\it cumulants},
$$
\tau_{l}(u|y)=\int_{\Scr}{[h(s)]}^{l}{\mbox e}^{-uh(s)}\rho(ds|y)
{\mbox { and }} \kappa_{l}(u)=\int_{\Scr}{[h(s)]}^{l}{\mbox
e}^{-uh(s)}{\rho}(ds)
$$
Define conditional distributions of the $(J_{j,n})$ given
$(U_n,\X)$, as \Eq\Pe(J_{j,n}\in ds|Y_{j},u)=\frac
{{[h(s)]}^{e_{j}}{\mbox
e}^{-us}\rho(ds|Y_{j})}{\tau_{e_{j}}(u|Y_{j})},{\mbox { for
}}j=1,\ldots, n(\p).\label{jumpd} \EndEq Note that, conditional on
$U_{n}$, each $J_{j,n}$ only depends upon $\X$ through
$(e_{j},Y_{j})$. In the homogeneous case their distributions no
longer depend upon $Y_{j}$ and correspond to the distribution of
$(J_{j,n})|U_{n},\p$ expressible as $\Pe(J_{j,n}\in ds|u)\propto
{[h(s)]}^{e_{j}}{\mbox e}^{-us}{\rho}(ds).$ We note further that
    $$
    \E\[h(J_{j,n})|u,\X\]=\frac{\tau_{1+e_{j,n}}(u|Y_{j})}{\tau_{e_{j,n}}(u|Y_{j})}
    \quad \mbox{and}\quad
    \E[h(J_{j,n})|u,\p]=\frac{\kappa_{1+e_{j,n}}(u)}{\kappa_{e_{j,n}}(u)}.
    $$
We now use these variables to describe the relevant posterior
distributions and related quantities.

\begin{thm}Let $P$ denote a NRM defined by the Poisson random measure
$N$
with mean intensity $\nu(ds,dx)=\rho(ds|x)\eta(dx)$. Let
$\X=(X_{1},\ldots, X_{n})$ denote a vector of random elements on a
Polish space such that $X_{1},\ldots,X_{n}|P$ are iid $P$, then
the following results hold.
\begin{enumerate}
\item[(i)]
The posterior distribution of $N$, given $(U_{n},\X)$, coincides
with the conditional distribution of the random measure
    $
    N^{*}_{n}=N_{n(\p)}+\sum_{j=1}^{n(\p)}\delta_{(J_{j,n},Y_{j})},
    $
where conditional on $U_{n}=u$ and $\X$, $N_{n(\p)}$ is a Poisson
random measure with intensity
    \begin{equation}
    \label{postlevy}
    \nu_{u}(\ddr s,\ddr x)={\mbox e}^{-uh(s)}\rho(\ddr s|x)\eta(\ddr
    x)={\mbox e}^{-uh(s)}F_{0}(dx|s)\rho(\ddr s),
    \end{equation}
not depending on $\X$, except through $U_{n}$.

\item[(ii)]Additionally, the $(J_{j,n})$ are, conditional on
$(U_n,\X)$, independent of $N_{n(\p)}$ and are mutually
independent with each $J_{j,n}$ having the distribution, specified
in~\mref{jumpd}

\item[(iii)] The posterior density of $U_{n}|\X$, is equivalent to
the density of
$\Gamma_{n}/[T_{n(\p)}+\sum_{j=1}^{n(\p)}h(J_{j,n})]|\X$, which is
$f_{U_{n}}(u|\X)\propto
u^{n-1}\edr^{-\psi(u)}\prod_{j=1}^{n(\p)}\tau_{e_{j}}(u|y_{j})$.
Similarly, the density of $U_{n}|\p$ is $f_{U_{n}}(u|\p)\propto
u^{n-1}\edr^{-\psi(u)}\prod_{j=1}^{n(\p)}\kappa_{e_{j}}
    (u).$
\end{enumerate}
    \end{thm}
Since $\mu$ and $P$ are functionals of $N$ the next two results
follow immediately.

\begin{prop}
The posterior distribution of $\mu|U_{n},\X$ is equivalent to the
conditional distribution, given $U_{n},\X$, of the random measure
$ \mu^{*}_{n}(dx)=\int_{0}^{\infty} sN^{*}_{n}(ds,dx)=
\mu_{n(\p)}(dx)+\sum_{j=1}^{n(\p)}h(J_{j,n})\delta_{Y_{j}} (dx), $
where conditional on $U_{n}$ and $\X$,
$\mu_{n(\p)}(dx):=\int_{0}^{\infty}sN_{n(\p)}(ds,dx)$ is a
completely random measure with L{\'e}vy measure specified in
\mref{postlevy}. This implies that the density of
$T_{n(\p)}|U_{n}=u,\X$ is $f_{T_{n(\p)}}(t|u)={\mbox
e}^{-ut}{\mbox e}^{\psi(u)}f_{T}(t)$.
\end{prop}

\begin{prop}
The posterior distribution of $P|\X$, is equivalent to the
conditional distribution of the random probability measure,
$$
P^{*}_{n}(dy)=
R_{n(\p)}P_{n(\p)}(dy)+\sum_{j=1}^{n(\p)}Q_{j,n}\delta_{Y_{j}}(dy),
$$
where $P_{n(\p)}(dx)=\mu_{n(\p)}(dx)/T_{n(\p)}$,
$R_{n(\p)}=T_{n(\p)}/T=1-\sum_{j=1}^{n(\p)}Q_{j,n}$ and
$Q_{j,n}=h(J_{j,n})/T$ for $j=1,\ldots,n(\p)$. The distribution of
all quantities, given $(U_{n},\X)$, is specified by Theorem 2.1.
\end{prop}
We now provide an initial description of $\Mcr$.
\begin{prop}
    The exchangeable marginal
distribution $\Mcr$ of the observations $\X$ can be represented as
    \begin{equation}
    \label{margin}
    \Mcr(\ddr x_1,\ldots,\ddr x_n)
    =\frac{1}{\Gamma(n)}\,
    \left\{\int_0^{+\infty}
    u^{n-1} \left[\prod_{j=1}^{n(\p)} \tau_{e_{j}}(u|y_j)\right]
    \,\edr^{-\psi(u)}\,\ddr u\right\} \,\prod_{j=1}^{n(\p)} \eta(\ddr
    y_j).
    \end{equation}
The corresponding EPPF is
    $$
    p(e_{1},\ldots,e_{n(\p)})={[\Gamma(n)]}^{-1}\int_0^{+\infty}
    u^{n-1} \left[\prod_{j=1}^{n(\p)} \kappa_{e_{j}}(u)\right]
    \,\edr^{-\psi(u)}\,\ddr u.
    $$
\end{prop}

\Remark In the case where $h(s)=s\in (0,\infty)$, the EPPF
appearing in Proposition 2.3 was first obtained by Pitman~(2003,
Corollary 6). Note that EPPF's appearing in Proposition 2.3 are
examples of {\it infinite} EPPF's. In the next section we shall
show that they may be represented as a mixture of tractable {\it
finite} EPPF's. The distinctions between finite and infinite
EPPF's are described in Pitman ~(2002, Section 2). \EndRemark

\section{Analysis of  the exchangeable marginal distribution $\Mcr$}
In this section we present a simpler description of the marginal
distribution $\Mcr$ and the corresponding EPPF. Such a description
facilitates its practical implementation and indeed has theoretical
implications as well. Note that ideally one would like an EPPF of
the form \Eq p(e_{1},\ldots,e_{n(\p)})
={v_{n,n(\p)}\prod_{j=1}^{n(\p)}w_{e_{j}}}. \label{Gibbs}\EndEq
where here $v_{n,n(\p)}$ is a positive quantity only depending on
$n$ and $n(\p)$, and the $w_{e_{j}}$ are a sequence of positive
numbers each only depending on $e_{j}$ for $j=1,\ldots, n(\p)$. One
aspect of such a representation is that it is easily sampled
according to variations of general Chinese restaurant processes.
Pitman~(2002) refers to such EPPF as having Gibbs form. However, it
is known [see Pitman~(2002), Theorem 42, p. 81] that the only
infinite EPPF admitting such a representation are the EPPF's derived
from a Dirichlet process and those derived from a Stable law of
index $0<\alpha < 1$. Among them, we mention the two-parameter
Dirichlet process and the generalized gamma class of processes. Here
again we show that the appropriate usage of the random variable
$U_{n}$, leads to tractable descriptions of the $\Mcr$ and the EPPF.

These simplifications are deduced from the following joint
distribution of  \\$(\p, T_{n(\p)},U_{n},\Y)$ given by $\Pe(\p,
T_{n(\p)}\in \ddr z,  U_{n}\in \ddr u, \Y \in \ddr \y)$ equal to,
\begin{equation}
\frac{1}{\Gamma(n)}\,u^{n-1}\edr^{-uz}f_{T}(z)\ddr u\ddr
z\prod_{j=1}^{n(\p)} \tau_{e_{j}}(u|y_j)\eta(\ddr
    y_j)
    \label{jointmarg}
\end{equation}
This distribution appears naturally in our derivation of the
posterior distribution given in the appendix. Now, for fixed $u>0$
and $\p$, let
    \Eq
    \PR(Y_{j}\in dy|U_{n}=u,\p):=H_{j,n}(\ddr y|u)=
    \frac{\tau_{e_{j,n}}(u|y)\eta(\ddr y)}{\kappa_{e_{j,n}}(u)}
    \qquad j=1,\ldots,n(\p)
    \label{margH}
    \EndEq

\begin{thm}
Let $\X$ denote the random variables with the exchangeable
distribution $\Mcr$ described in Proposition 2.3. Then the
distribution of $\X$ may be described as follows.
\begin{enumerate}
\item[(i)] The distribution of $\X|U_{n}=u,\p$ is such that the unique values
$Y_{1},\ldots,Y_{n(\p)}$ are independent with respective
distributions given by~\mref{margH}. In the homogeneous case, it
follows that $Y_{1},\ldots,Y_{n(\p)}|\p$ are iid with common
distribution $\eta(dx)/\eta(\Xcr)=H(dx)$ not depending on $U_{n}$.
\item[(ii)]The distribution of $\p|U_{n}=u$,
    $\Pe(\p=\{C_{1},\ldots,C_{n(\p)}\}|U_{n}=u)$, is
    a conditional finite EPPF given by
    \begin{equation}
    p(e_{1},\ldots,e_{n(\p)}|u)=\frac{
    \edr^{-\psi(u)}\prod_{j=1}^{n(\p)}\kappa_{e_{j}}
    (u)}{\int_{0}^{\infty}t^{n}\edr^{-ut}f_{T}(t)dt}.
    \label{part}
    \end{equation}
where
$E[T^{n}_{n(\p)}|U_{n}=u]=\edr^{\psi(u)}{\int_{0}^{\infty}t^{n}\edr^{-ut}f_{T}(t)dt}$
That is, conditional on $U_{n}=u$, $\p$ is a finite Gibbs partition.
\item[(iii)] The marginal density of $U_{n}=\Gamma_{n}/T$ is
    $$
    f_{U_{n}}(u)={[\Gamma(n)]}^{-1}u^{n-1}\int_{0}^{\infty}t^{n}\edr^{-ut}f_{T}(t)dt.
    $$
\end{enumerate}
\end{thm}
\Proof Statement (i) follows by applying Bayes rule
to~\mref{jointmarg}. Statement (iii) is straightforward. An
application of Bayes rule also yields readily a description of the
distribution of $\p|U_{n}=u$, where the normalizer is $\sum_{\p}
\prod_{j=1}^{n(\p)}\kappa_{e_{j}}(u)$. Here $\sum_{\p}$ stands for
the sum over all partitions of the set of integers $\{1,\ldots,n\}$.
The simpler form in~\mref{part} may be obtained by noting some known
relationships between cumulants, partitions and moments. However,
for immediate clarity one can use~\mref{jointmarg} to establish the
identity $f_{U_{n}}(u)={[\Gamma(n)]}^{-1}u^{n-1}{\mbox
e}^{-\psi(u)}\sum_{\p} \prod_{j=1}^{n(\p)}\kappa_{e_{j}}(u)$. The
result then follows by noting the form of $f_{U_{n}}(u)$ given in
(iii). \EndProof

The next proposition offers another description of the
distribution of the unique values.

\begin{prop} Suppose that $\X$ has distribution $\Mcr$. Then using the decomposition
$\nu(ds,dx)=F_{0}(dx|s)\rho(ds)$, it follows that the distribution
of $Y_{1},\ldots,Y_{n(\p)}$ given
$U_{n},J_{1,n},\ldots,J_{n(\p),n},\p$ is such that the $Y_{j}$ are
independent with respective distributions $\Pe(Y_{j}\in
dy|J_{j,n})=F_{0}(dy|J_{j,n})$. That is, the conditional
distributions only depend on $U_{n}$ through the $(J_{j,n})$
\end{prop}

\Remark The above results demonstrate the important role that
$U_{n}$ plays in simplifying the above quantities. In effect,
conditioning on $U_{n}$, reveals {\it conditional likelihoods}
that have exponential form. This exponential form bears
resemblance to those appearing in the posterior analysis of
neutral to the right processes and the Levy moving average models
discussed in James~(2005a). Models of this type are most amenable
to the direct application of James~(2005a, Proposition 2.3). It is
then not surprising that one may notice some similarities between
our posterior characterizations and those described by
Doskum~(1974), Ferguson~(1974), Hjort(1990) and Kim~(1999) for NTR
models. We point out however that these models are otherwise very
different, see James~(2005a,2005c). \EndRemark

\Remark We note further that conditioning on $T$, instead of
$U_{n}$, does not lead to simplified expressions. See however
Pitman~(2003) for important interpretations of conditioning on $T$
in the stable case. \EndRemark
 We close this section by showing how $U_{n}$ is related to
what is called the {\it real inversion formula for the Laplace
transform} as described in Bondesson~(1992, eq. (6.2.1), p.
92)[see Feller~(1971,VII)]. In some sense it shows that $U_{n}$ is
asymptotically sufficient for $T$. The result below follows
directly from Feller~(1971, VII).

\begin{prop} Let $Y_{n}=n/U_{n}=nT/\Gamma_{n}$. Then the pdf of $Y_{n}$ is given by
\Eq
f_{Y_{n}}(y)=\frac{1}{\Gamma(n+1)}{\(\frac{n}{y}\)}^{n+1}\int_{0}^{\infty}{\mbox
e}^{-nt/y}t^{n}f_{T}(t)dt=\frac{{(-1)}^{n}}{n!}{\(\frac{n}{y}\)}^{n+1}\phi^{(n)}(n/y),
\label{inverse} \EndEq where $\phi^{(n)}$ denotes the $n$-th
derivative of $\phi$. Additionally it follows that $f_{Y_{n}}(t)$
converges to $f_{T}(t)$ uniformly in every finite interval, as
$n\rightarrow \infty.$ Hence ~\mref{inverse} is an inversion
formula for the Laplace transform of $T$
\end{prop}
\Remark Comparing~\mref{part} with~\mref{Gibbs} shows that
conditional on $U_{n}$, the distribution of $\p$, is a finite EPPF
of  Gibbs form. In particular for fixed $u$, we see that in this
case,
$$
v_{n,n(\p)}=
E[T^{n}_{n(\p)}|U_{n}=u]=\edr^{\psi(u)}{\int_{0}^{\infty}t^{n}\edr^{-ut}f_{T}(t)dt},$$
which importantly does not depend on $n(\p)$, and
$w_{e_{j}}=\kappa_{e_{j,n}}(u)$. This has many interesting
consequences, of which we shall highlight a few in the forthcoming
sections. \EndRemark

\subsection{Distributional results and moment formulae for complex functionals}
The next result gives an expression for the distribution of
$n(\p)$ given $U_{n}$.
\begin{prop}
Let $n_{1}\ge n_{2}\ldots \ge n_{k}\ge 1$ denote an ordering
(composition) of the corresponding $(e_{1},\ldots, e_{k})$. Recall
that $n(\p)$ represents the number of unique values of
$(X_{1},\ldots,X_{n})$. Then under the distribution $\Mcr$, the
conditional distribution  of $n(\p)$ given $U_{n}=u$, is given by
$$
\PR(n(\p)=k|U_{n}=u)=\frac{
\edr^{-\psi(u)}}{\int_{0}^{\infty}t^{n}\edr^{-ut}f_{T}(t)dt}
\frac{n!}{k!}\sum_{(n_{1},\ldots,n_{k})}
\prod_{j=1}^{k}\frac{\kappa_{n_{j}}(u)}{n_{j}!}
$$
for $k=1,\ldots, n$. Where the sum corresponds to the sum over all
compositions $(n_{1},\ldots,n_{k})$ of $n$ of size $n(\p)=k$.
\end{prop}

The fact that $\p|U_{n}$ is a finite Gibbs partition allows us to
apply Pitman~(2002, eq. 98) [see also Kolchin~(1986)] to immediately
deduce the following generalization of the Ewens sampling formula,
Ewens (1972), and equivalently Antoniak~(1974, Proposition~3).
\begin{prop}
Define a random vector $(|\Pi_{i,n}|, for 1\leq i\leq n)$ of
non-negative counts by $|\Pi_{i,n}|=\sum_{j=1}^{n}I(|C_{j,n}|=i)$
for $i=1,\ldots, n$. Where specifically, $|\Pi_{i,n}|$ denotes the
number of cells of size $i$. The distribution of $(|\Pi_{i,n}|,
{\mbox { for }} 1\leq i\leq n)$ given $U_{n}=u$, can be
represented as \Eq \Pe(|\Pi_{j,n}|=m_{j}, 1\leq j\leq
n|u)=\frac{n!
\edr^{-\psi(u)}}{\int_{0}^{\infty}t^{n}\edr^{-ut}f_{T}(t)dt}\prod_{j=1}^{n}{\left(\frac{\kappa_{j}
    (u)}{j!}\right)}^{m_{j}}\frac{1}{m_{j}!}
\label{EwensAnt}
    \EndEq
where $\sum_{j=1}^{n}m_{j}=k$ and $\sum_{j=1}^{n}jm_{j}=n.$
Equivalently~\mref{EwensAnt} is the conditional distribution, given
$U_{n}$, of the number of values of $(X_{1},\ldots,X_{n})$ appearing
$1$ time, $2$ times etc, corresponding to the numbers
$(m_{1},\ldots, m_{n})$, when $\X$ has distribution $\Mcr$.
\end{prop}

The recognition that $\Mcr$ is the $n-{\mbox {th}}$ moment measure
of $P$ allows one to obtain easily otherwise complex expressions
for moments of functionals of $P$. The discussion in Ishwaran and
James (2003a, section 3.2) combined with Proposition 2.1 leads to
the following formula. Recall that $\sum_{\p}$ denotes the sum
over all partitions of the set $\{1,\ldots,n\}$.

\begin{prop}
Suppose that $g_{1},\ldots, g_{n}$ are real--valued functions on
$\X$ and define the functionals $P(g_{i})=\int_{\mathcal
X}g_{i}(x)P(dx)$ for $i=1,\ldots n$. Assume that for each $n,$
$E\left[\prod_{i=1}^{n}{(P(g_{l}))}\right]< +\infty$. \Enumerate
\item[(i)] Then, $\E\left[\prod_{l=1}^{n}{P(g_{l})}\right]$ is
equal to,
\begin{equation}
    \label{eq:mixture}
\int_{0}^{\infty}\left[\sum_{\p}\pi(\p|u)
    \prod_{j=1}^{n(\p)}\int_{\X}\[\prod_{i \in C_{j}}
    {g_{i}(y)}\]H_{j,n}(\ddr y|u)\right]f_{U_{n}}(u)\,\ddr
    u.
    \end{equation}
    \item[(ii)]For integers $n_{1},\ldots, n_{q}$ chosen such that
$\sum_{i=1}^{q}n_{i}=n$ for an integer $q\leq n$, it follows that,
$E\left[\prod_{l=1}^{q}{(P(g_{l}))}^{n_{l}}\right]$ is equal to
    \begin{equation}
    \label{eq:mixture}
   \int_{0}^{\infty}\left[\sum_{\p}\pi(\p|u)
    \prod_{j=1}^{n(\p)}\int_{\X}\[\prod_{l=1}^{q}
    {[g_{l}(y)]}^{e_{j,l}}\]H_{j,n}(\ddr y|u)\right]f_{U_{n}}(u)\,\ddr
    u.
    \end{equation}
    Where $e_{j,l}$ denotes the number of indices in $C_{j}$
    associated with $g_{l}$.
    \EndEnumerate
\end{prop}

\Remark It is interesting to note that all our results conditioned
on $U_{n}$, contain the known unconditional results for the
Dirichlet process. This is because the Dirichlet process is
independent of $U_{n}$. To see this, note that the Dirichlet
process with total mass $\theta$, corresponds to the choice of
$\rho(ds)=\theta s^{-1}{\mbox e}^{-s}ds$. It follows that for each
$j$ that $\kappa_{j}(u)=\theta (1+u)^{-j}\Gamma(j)$ and
$E[T_{n(\p)}|u]=(1+u)^{-n}[\Gamma(\theta)/\Gamma(\theta+n)]$.
Aditionally the $f_{U_{n}}(u|\p):=f_{U_{n}}(u)\propto
u^{n-1}(1+u)^{-(n+\theta)}$. That is $U_{n}=\Gamma_{n}/T$ is a
gamma-gamma random variable independent of $\X$. Equivalently
$1/(1+U_{n})$ is a ${\mbox {Beta}}(\theta,n)$ random variable.
Hence,~\mref{EwensAnt}, specializes to
$$
\Pe(|\Pi_{j,n}|=m_{j}, 1\leq j\leq n|u)=\frac{n!
}{\prod_{i=1}^{n}(\theta+i-1)}\prod_{j=1}^{n}{\left(\frac{\theta
    }{j}\right)}^{m_{j}}\frac{1}{m_{j}!}.
$$

This equates to the Ewens sampling formula derived by Ewens
(1972), which is equivalent to the result in Antoniak (1974,
Proposition 3) describing the number of values of
$(X_{1},\ldots,X_{n})$ appearing $1$ time, $2$ times etc,
corresponding to the numbers $(m_{1},\ldots, m_{n})$.
Additionally, note that~\mref{part} becomes, \Eq {\mbox
{PD}}(\p|\theta):=\frac{\theta^{n(\p)}\prod_{j=1}^{n(\p)}(e_{j}-1)!}
{\prod_{i=1}^{n}(\theta+i-1)}. \label{CREPPF}\EndEq which is the
variant of Ewens sampling formula, often called the Chinese
restaurant process. [See Pitman (2002, p. 60) and Ishwaran and
James~(2003a)]. The calculations for the Dirichlet process
involving $U_{n}$ may be found in James~(2005b), where it is shown
that $U_{n}$ and its variants still play a significant role.
\EndRemark

\section{Mixture models} In terms of statistical applications,
owing to the success of the Dirichlet process, one of the most
fruitful ways for exploiting NRM's is their potential use as basic
building blocks in hierarchical mixture models. In this setting,
$\Xn$ are missing values which capture the clustering structure
within the data. This class of models was first introduced, for
the Dirichlet process, by Lo (1984) and later popularized by the
development of suitable MCMC techniques in Escobar and West
(1995). See  Dey, M\"uller and Sinha (1998) and Ishwaran and James
(2001, 2003a,b) for subsequent developments. Recently, mixtures of
Dirichlet process have been generalized to mixtures of
stick--breaking priors in Ishwaran and James (2001, 2003b) and
random measures driven by increasing additive processes in
Nieto--Barajas, Pr\"unster and Walker (2004). A recent example of
application of this class is provided in Lijoi, Mena and
Pr\"unster (2004) where the clustering behaviour is modeled
according to a normalized inverse Gaussian process. Ishwaran and
James~(2003a) also introduce a general class of species sampling
mixture models and describe various algorithms for efficient
implementation. See also Hoff (2003, section 4) for an interesting
use of the Dirichlet process mixture model framework. Those ideas
naturally extend to models based on the NRM's.

We first recall the model as set up by Lo (1984). Suppose
$\{f(\,\cdot\,|x):\, x\in\Xcr\}$ is a family of non--negative
kernels defined on a Polish space ${\mathbf W}$ such that
$\int_{{\mathbf W}}f(w|x)\,\lambda(\ddr w)=1$ for any $x$ in
$\Xcr$ and for some $\sigma$--finite measure $\lambda$. Next, let
$\W=(W_{1},\ldots, W_{n})$ be a vector of $\Wbb$--valued random
elements such that, given $\Xn$ from a NRMI $P$, they are
independent and $W_j$ admits density, with respect to $\lambda$,
$f(\,\cdot\,|X_j)$. This is the same as supposing that
$W_1,\ldots,W_n$ are exchangeable draws from the random density
$\tilde f(\,\cdot\,)=\int_\X f(\,\cdot\,|x)\,P(\ddr x)$. One is
naturally interested in the determination of the distribution of
the posterior density $\tilde f$, given the observations $\W$.
However, one gains more flexibility in working directly with the
posterior distribution of $P$ or $\mu$ given $\W$. That is,
$\tilde {f}$, is then seen as one of many possibly interesting
functionals of $P$. Moreover, under certain identifiability
assumptions, the estimation of the mixing distribution $P$ is of
primary concern.

Notice that the above description shows that the joint
distribution of $(\W,\X,P,\U_{n})$ can be written as,
$$
\[\prod_{i=1}^{n}f(W_{i}|X_{i})\]\PR(dP|\X)\Mcr(d\X|u)f_{U_{n}}(u),
$$
where $\PR(dP|\X)$ denotes the posterior distribution of $P$
described by Theorem 2.1 or Proposition 2.2. One could then apply
arguments similar to those exploited in Lo~(1984) and in Ishwaran
and James~(2003a) to yield analogous characterizations of the
posterior distribution. We shall not present those here. In the
next section we will describe a general Monte-Carlo which can be
used to sample from the posterior distribution. For a better
understanding of this connection we note that an application of
Proposition 2.4 shows that the marginal distribution of $\W$ is
given by
    $$
    \int_{0}^{\infty}\left[\sum_{\p}\pi_n(\p|u)
    \prod_{j=1}^{n(\p)}\int_{\Xcr}\left[\prod_{i\in
    C_{j}}f(W_{i}|y)\right]\, H_{j,n}(\ddr y|u) \right]\, f_{U_{n}}(u)\,\ddr u.
    $$
This is a special case of~\eqref{eq:mixture}, with
$g_{i}(y)=f(W_{i}|y)$.

\Remark It is important to note that in mixture models both $\X$
and $U_{n}$ are viewed as missing values. Hence it is quite
natural to work with the distribution of $\X|U_{n}$ as a prior and
subsequently $\X|U_{n},\W$ \EndRemark

\section{Sampling from $\Mcr$ and related functionals}
As mentioned earlier, obtaining a tractable form of the marginal
distribution $\Mcr$ is crucial to both practical implementation and
theoretical understanding of these models. In particular,
understanding how to sample $\Xn$ from $\Mcr$ is important for
applications involving mixture models. We discuss briefly some ideas
on how this may be done. One aspect of our expressions is the
appearance of the cumulants $\kappa_{l}(u)$ and the corresponding
moments,
    $$
    m_{n}(u)=E[T^{n}_{n(\p)}|U_{n}=u].
    $$
In many cases either the cumulants are easy to calculate or the
moments are. Moreover, one can use the following result of Theile in
order to recover one from the other. That is, for any integer $n$
    \Eq
    \label{theile}
    \kappa_{n}(u)=m_{n}(u)-\sum_{l=1}^{n-1}\binom{n-1}{l-1}
    \kappa_{l}(u)m_{n-l}(u).
    \EndEq
Many mathematical packages can easily deal with \mref{theile}.
Similar to the case of the Dirichlet process, it is noted that
many complex expressions can be approximated by obtaining draws
from $\Mcr$. Using Proposition~2 and Corollary~1 a draw from
$\Mcr$ may be conducted as follows. First one draws
$U_{n}=\Gamma_{n}/T$, either directly or by drawing from the
independent pair $(\Gamma_{n},T)$ according to the gamma density
of $\Gamma_{n}$ and the density $f_{T}$. Given $U_n=u$, one draws
$\p$ from $\pi_n(\p|u)$. Given $\p$ and $U_n=u$, one finally draws
$\X$ which amounts to independently sampling $Y_{j}$ from
$H_{j,n}(\,\cdot\,|u))$, for $j=1,\ldots,n(\p)$. Since the
normalizing constant, $m_{n}(u)$, in $\pi(\p|u)$ is fairly simple,
one may often be able to devise a simple scheme to draw $\p$,
given $U_n=u$, exactly. If this is not the case, one can use a
simple variation of a weighted Chinese restaurant (WCR) process
[see Lo, Brunner, and Chan~(1996) and Ishwaran and James~(2003a)],
which can be deduced from~James (2002, Lemma~2.3).

\subsection{Generalized Chinese Restaurant and P\'olya Urn procedures}
Here, we use the fact that these models are structurally similar
to those discussed in section 4 of James~(2005a). It follows that
one can use section 4.4. of James~(2005a) to deduce general
extensions of P\'olya Urn Gibbs samplers and SIS procedures given
by Escobar(1994), Liu~(1996), and the Gibbs sampling/SIS
procedures based on a generalized Chinese restaurant process as
mentioned above. As such, we will only sketch out the relevant
probabilities and refer the reader to James~(2005a), and
references therein, for additional mechanics of the
implementation.

For greatest flexibility we will give the relevant probabilities
needed to sample approximately from models related to a joint
density proportional to
$$
\[\prod_{i=1}^{n}g_{i}(X_{i})\]\Mcr(d\X|u)
$$
which is deducible from Proposition 3.5. We note that
$\Mcr(d\X|u)$ has an urn type representation which can be deduced
from an application of James (2005a, Proposition 5.1). Similar to
James~(2005a, equation 40), define for $r=0,\ldots, n-1$
conditional probabilities,
    $$
    \Pe(X_{r+1}\in
    dx|\X_{r},u)=\frac{l_{0,r}}{c_{r}}\lambda_{r}(dy|u)+\sum_{j=1}^{n(\p_{r})}
    \frac{l_{j,r}(Y_{j})}{c_{r}}\delta_{Y_{j}}(dx)
    $$
where $\X_{r}=\{X_{1},\ldots,X_{r}\}$, $\lambda_{r}(dx|u)\propto
g_{r+1}(x)H_{1,1}(dx|u)$,
    $$
    l_{0,r}=\int_{\Xcr}g_{r+1}(y)\, H_{1,1}(\ddr
    y|u)\quad\mbox{and}\quad
    l_{j,r}(y)=g_{r+1}(y)\tau_{1+e_{j,r}}(u|y)/\tau_{e_{j,r}}(u|y).
    $$
Additionally $c_{r}=l_{0,r}+\sum_{j=1}^{n(\p_{r})}l_{j,r}(Y_{j})$.
Examining James~(2005a, section 4.4.) we see these are the
ingredients to implement general analogues of the P\'olya Urn
Gibbs Sampler and SIS procedures described by Escobar~(1994) and
Liu~(1996). To get the Chinese restaurant type procedures one
samples partitions $\p$ based on probabilities derived from
$l_{0,r}$ and
$$
l_{j,r}=\int_{\Xcr}l_{j,r}(y)\,
    \prod_{i\in C_{j,r}}g_{i}(y)H_{j,r}(\ddr y|u){\mbox { for }}
    j=1,\ldots,n(\p_{r})
$$
where $\p_{r}$ denotes a partion of the integers $\{1,\ldots, r\}$
and each $C_{j,r}=\{i\leq r: X_{i}=Y_{j}\}$ denotes the
corresponding cells. Additionally,
$l(r)=l_{0,r}+\sum_{j=1}^{n(\p_{r})}l_{j,r}$. In particular,
applying the SIS WCR procedure described in James~(2005a), now
leads to a sampling $\p$ from a density $q(\p|u)$, which satisfies
$$
L(\p|u)q({\p}|u)=\pi(\p|u)\prod_{j=1}^{n(\p)}\int_{\Xcr}\prod_{i\in
C_{j}}g_{i}(y)H_{j:n}(dy|u)
$$
where $L(\p|u)=\prod_{r=1}^{n}l(r-1)/m_{n}(u)$. This is justified
by James~(2002, Lemma 2.3). Note that setting
$g_{i}(X_{i})=f(W_{i}|X_{i})$ leads to sampling procedures for
mixture models. Setting $g_{i}(x)=1$ leads to sampling from
$\Mcr$. In particular this algorithm includes the classical
Chinese restaurant process for the Dirichlet process.

\section{Illustrative Examples} Here we study two examples which are
connected to the Dirichlet process, but require a more delicate
analysis. In section 7, we address a more involved class of models.
Hereafter, we let $\Bcr(a,b)$, denote the fact that a random
variable has a Beta distribution with parameters $(a,b)$. Let
$\Bcr(x|a,b)$ denote its density. Similarly $\Gcr(a)$ denotes the
law of a gamma random variable with shape $a$ and scale $1$. $G_{a}$
denotes the corresponding gamma random variable having density
$f_{G_{a}}(x)=\Gcr(x|a)$.

\medskip
\subsection{Classes of Dependent Dirichlet processes}
Here we present a large class of models which share the same EPPF
as the Dirichlet process but are otherwise substantially more
complex. This class is seen to add more flexibility to the
Dirichlet process, and may be of particular interest in mixture
models. One can also see that a study of a subclass of such models
is also related to the exposition in Aldous and Pitman~(2002).

That is we build NRM based on \Eq
\nu_{0,\theta}(ds,dx)=F_{0}(dx|s)\theta s^{-1}{\mbox
e}^{-s}.\label{DPint}\EndEq It is evident that $T\overset
{d}=G_{\theta}$. A careful examination of~\mref{DPint} shows that a
class of dependent Dirichlet processes may be described in terms of
a stick-breaking representation,
    \Eq
    P_{\theta}(\cdot)=\sum_{i=1}^{\infty}{\tilde
    Q}_{i}\delta_{Z_{i}}(\cdot)\qquad {\mbox { with }} \qquad {\tilde
    Q}_{i}=W_{i}\prod_{j=1}^{i-1}(1-W_{j}). \label{DPrep}
    \EndEq
The $W_{i}$ are independent ${\mbox {Beta}}(1,\theta)$,
corresponding to the usual representation of Sethuraman~(1994), but
now the $(Z_{i})$ are no longer independent of the $(W_{i})$. A
technical point is that $\tilde
Q_{i}=W_{i}\prod_{j=1}^{i-1}(1-W_{j})$ are the points ranked by
size-biased sampling and the corresponding $Z_{i}$ has conditional
distribution depending on $T{\tilde Q}_{i}={\tilde J_{i}}$, i.e.
$\Pe(Z_{i}\in dz|{\tilde J}_{1},{\tilde J}_{2},\ldots
)=F_{0}(dz|{\tilde J}_{i})$. Note importantly that the distribution
of ${\tilde J}_{i}$ is much more manageable than the distribution of
the ranked points of a gamma process.

We can describe the posterior distribution in the following way.
First, notice that
    \Eq
    \nu_{u}(ds,dx)={\mbox
    e}^{-us}F_{0}(dx|s)\theta s^{-1}{\mbox e}^{-s}
    \label{postLDP},
    \EndEq
which implies that $(1+U_{n})T_{n(\p)}\overset {d}=G_{\theta}$  and
$ (1+U_{n})\mu_{n(\p)}(dx)\overset {d}=G_{\theta}\sum_{i=1}^{\infty}
{\tilde Q}_{i}\delta_{Z_{i,n}}$. Here $(\tilde {Q}_{i})$ have the
same distribution as described above in~\mref{DPrep}, they are
independent of $\X$ and $U_{n}$, while $(Z_{i,n})$ are now random
variables depending on $U_{n}$. Specifically, conditioned on the
sequence $(U_{n},G_{\theta}{\tilde Q}_{1},G_{\theta}{\tilde
Q}_{2},\ldots)$, the $(Z_{i,n})$ are independent with distributions
    \Eq \PR(Z_{i,n}\in dx|U_{n},G_{\theta}{\tilde
    Q}_{1},G_{\theta}{\tilde Q}_{2},\ldots)= F_{0}(dx|G_{\theta}{\tilde
    Q}_{i}/(1+U_{n}))
    \label{distZ}
    \EndEq
where $G_{\theta},U_{n}$ and $({\tilde Q}_{i})$ are independent. Now
setting $(1+U_{n})J_{j,n}=G_{j,n}$ and one has the conditional
distributions of $G_{j,n}|U_{n},\X$ and $G_{j,n}|U_{n},\p$ as,
    \Eq
    \PR(G_{j,n}\in ds|Y_{j},u)\propto
    F_{0}(dY_{j}|s/(1+u))\Gcr(s|e_{j,n})ds
    \label{distG}
    \EndEq
and
    $$
    \PR(G_{j,n}\in ds|u)=\Gcr(s|e_{j,n})ds.
    $$
Additionally, set $T_{n}=(1+U_{n})T$. That is $T_{n}\overset
{d}=G_{\theta}+\sum_{j=1}^{n(\p)}G_{j,n}$. The conditional density
of $V_{n}:=1/(1+U_{n})$, given $\X$ is specified by
    \Eq
    f_{V_{n}}(v|\X)\propto
    v^{\theta-1}{(1-v)}^{n-1}\prod_{j=1}^{n(\p)}
    \int_{0}^{\infty}F_{0}(dY_{j}|sv)\Gcr(s|e_{j,n})ds.
    \label{vv}
    \EndEq
$V_{n}|\p$ is a $\Bcr(\theta,n)$ random variable, independent of
$\p$. We now summarize some facts about this process in the next
result.

\begin{thm} Suppose that $P_{\theta}$ denotes a class of dependent
Dirichlet processes defined as in ~\mref{DPrep} via the
intensity~\mref{DPint}. Note that $V_{n}:=1/(1+U_{n})$ given $\X$
has distribution~\mref{vv}. Then the following results hold
\Enumerate
\item[(i)]The
posterior distribution of $P_{\theta}|V_{n},\X$ is equivalent to
the random probability measure
$$
\frac{G_{\theta}}{T_{n}}\sum_{i=1}^{\infty}{\tilde
Q}_{i}\delta_{Z_{i,n}}+ \sum_{j=1}^{n(\p)}Q_{j,n}\delta_{Y_{j}}
$$
where the $({\tilde Q}_{i})$ are equivalent in distribution to
those in ~\mref{DPrep}, and are independent of $\X$. The sequence
$(Z_{i,n})$ has distribution specified by~\mref{distZ}.
Additionally $Q_{j,n}=G_{j,n}/T_{n}$, where the distribution of
$G_{j,n},T_{n}|V_{n},\X$ is specified by~\mref{postLDP}
and~\mref{distG}.
\item[(ii)]The distributions of $((G_{j,n}),T_{n},(Q_{j,n}))|\p$, are
the same as given $V_{n}$ and $\p$, and equate to the classical
results for the Dirichlet process posterior distribution. That is,
given $\p$, $(Q_{1,n},\ldots,Q_{n(\p),n})$ is a Dirichlet
$(e_{1},\ldots,e_{n(\p)};\theta)$ vector. Equivalently, each
$Q_{j,n}$ is $\Bcr(e_{j,n},\theta+n-e_{j,n})$, since the
$(G_{j,n})$ given $\p$ are independent $\Gcr(e_{j,n},1)$ and
$T_{n}|\p$ is $\Gcr(\theta+n)$. \EndEnumerate
\end{thm}
It is evident again that the prediction rule does not in general
have a simple form. However, the next result yields a nice
description of the marginal distribution of $\X$.
\begin{prop} Suppose that $P_{\theta}$ is a dependent Dirichlet
process as described above. Then the marginal distribution of
$\X|V_{n}=v$, where $V_{n}$ is $\Bcr(\theta,n)$, is given by \Eq
PD(\p|\theta)\prod_{j=1}^{n(\p)}\int_{0}^{\infty}F_{0}(dy_{j}|sv)\Gcr(s|e_{j,n})ds
\label{jointDDP}\EndEq The expression in~\mref{jointDDP} shows
that the conditional distribution of
$Y_{1},\ldots,Y_{n(\p)}|V_{n}=v,\p$ are independent with
respective distributions
    $$
    \Pe(Y_{j}\in
    dy|v)=\int_{0}^{\infty}F_{0}(dy_{j}|sv)\Gcr(s|e_{j,n})ds.
    $$
Additionally, $Y_{1},\ldots,Y_{n(\p)}|(G_{j,n}),V_{n}=v,\p$ are
independent with distributions $F_{0}(dy_{j}|G_{j,n}v)$.
\end{prop}

Notice that in every case the distribution of $(Y_{j})|V_{n},\p$
are expressible via gamma mixing measures. The next result
demonstrates a particularly simple case.
\begin{prop} Suppose that for $a>0, \delta>0$,
$F_{0}(dy|a)\propto a^{\delta}y^{\delta-1}{\mbox e}^{-ay}dy$
corresponds to a gamma random variable. Then
$Y_{1},\ldots,Y_{n(\p)}|V_{n},\p$ are independent with each
$Y_{j}\overset {d}=R_{j}/V_{n}$, where $R_{j}$ is independent of
$V_{n}$ and the distribution of $R_{j}|\p$ is
$$
f_{R_{j}}(r|\p)\propto
r^{\delta-1}{(1+r)}^{-(\delta+e_{j,n})}{\mbox { for }}0<r<\infty
$$
Or equivalently $1/(1+R_{j})$ given $\p$ is
$\Bcr(e_{j,n},\delta)$. Additionally the $Z_{i,n}$ appearing in
Proposition 6.1 satisfy $Z_{i,n}\overset {d}={\tilde
R}_{i}/({\tilde Q}_{i}V_{n})$ where $1/(1+{\tilde {R}}_{i})$ is
independent of $({\tilde Q}_{i}), V_{n}$ and has a
$\Bcr(1,\delta)$ distribution. The distribution of the $Z_{i}$ is
given by setting $V_{0}:=1$
\end{prop}
\Remark Note that unlike the case of the usual Dirichlet process
$V_{n}$ is not independent of $\X$, however its distribution is
independent of $\p$ and the marginal distribution of $V_{n}$ is
the same in both the dependent and classical case. Proposition 6.2
shows that it is easy to sample from $\Mcr$. Since $V_{n},\p$ are
independent, one can draw from $PD(\p|\theta)$ according to the
classical Chinese restaurant, then sample independently $V_{n}$
from a $\Bcr(\theta,n)$. It then remains to draw the
$Y_{1},\ldots,Y_{n(\p)}$ conditional on $V_{n},\p$. Or as an
intermediate step, conditional on $(G_{j,n}),V_{n}=v,\p$.
\EndRemark

\subsection{NRMs derived from Beta processes} Suppose that $\rho(\ddr
s|x)=s^{-1}(1-s)^{c(x)-1} ds\, c(x)$, for $0<s\leq 1$. The
corresponding process $\mu$ is a \textit{beta process}. Exploiting
beta processes as an ingredient for constructing NRM models leads
to a different posterior behaviour from the one analyzed in Hjort
(1990). The L\'evy measure associated with $\mu_{u}$ in
\mref{postlevy} is
    $$
    \nu_{u}(\ddr s,\, \ddr x)=\edr^{-us}\,
    s^{-1}(1-s)^{c(x)-1}\ddr s\, \,\eta(\ddr x)
    \qquad {\mbox
    { for }}\quad s\in(0,1]
    $$
and, hence, $\mu_{u}$ is not a beta process. Additionally, the
distribution of the jumps $J_{j,n}$, given $(U_{n},\X)$ is
    $$
    \Pe(J_{j,n}\in\ddr s|Y_{j},u)=\frac{1}{_{1}F_{1}(e_{j,n},c(y_{j})+e_{j,n};-u)}
    \, {\mbox
    e}^{-us}\, s^{e_{j,n}-1}\, (1-s)^{c({y_{j}})-1}\, \ddr s
    $$
where $_1F_1$ is the {\it confluent hypergeometric function}. It
follows that one can write $m(\ddr x_{i}|x_{1},\dots, x_{i-1},u)$
as,
    \begin{multline*}
    _1F_1(1,c(x_i)+1;-u) \, \eta(\ddr x_i)+\\
    \sum_{j=1}^{n(\p_{i-1})}\frac{e_{j,i-1}}{e_{j,i-1}+c(y_j)}
    \frac{_{1}F_{1}(e_{j,i-1}+1,c(y_j)+e_{j,i-1}+1;-u)} {
    _{1}F_{1}(e_{j,i-1},c(y_j)+e_{j,i-1};-u)}\, \delta_{y_{j}}(\ddr
    x_i).
    \end{multline*}
Finally, the marginal distribution of $\X$, given $U_n$, can be
represented as
    $$
    {[m_{n}(u)]}^{-1}\left[\prod_{j=1}^{n(\p)}\frac{\Gamma(e_{j,n})
    \Gamma(c(y_j)+1)}{\Gamma(e_{j,n}+c(y_j))}\right]
    \prod_{j=1}^{n(\p)}
    {_{1}F_{1}(e_{j,n},c(y_j)+e_{j};-u)}\,\eta(\ddr y_{j}).
    $$
where $m_{n}(u)$ can be expressed using \mref{theile}, where
$\kappa_{1}(u)=m_{1}(u)$.

\Remark We point that by setting $c(s):=1$ and letting
$\eta(dx)=\theta H(dx)$, yields models defined by the scale
invariant Poisson process, which is of importance in a variety of
applications. In particular it is known that $T$ has the important
Dickman distribution. Moreover, in this case, the class of NRM have
been discussed previously. Specifically, Arratia, Barbour and
Tavare~(1999, 2004) show that the distribution of $P|T\leq 1$ is in
fact a Dirichlet process with shape $\theta$. See Arratia, Barbour
and Tavare~(2004) for further implications and details. \EndRemark

\section{Generalized Gamma Convolution processes: NRM related to Dirichlet mean
functionals} This last example demonstrates the great flexibility of
the $h$-biased framework, which allows us to describe a large class
of NRM in terms of the tractable $({\tilde {Q}}_{i})$ of the
Dirichlet process and related variables. These models will be based
on random variables $T$ which have distributions that are {\it
Generalized Gamma Convolutions} (GGC).  This large class of {\it
self-decomposable} infinite-divisible random variables was
introduced by Thorin~(1977, 1978), further developments are given in
Bondesson~(1979, 1992). In particular, a subclass of such models
will be seen to be connected to the study of mean functionals of a
Dirichlet process, initiated by Cifarelli and Regazzini~(1990). We
believe that our discussion will also shed some new light on these
two lines of research which are essentially duals of one another.
Moreover, our approach shows that one may extend the study of mean
functionals of Dirichlet processes to this larger and more flexible
setting. Quite strikingly, this class of NRM is very rich, including
the stable law processes of index $0<\alpha<1$, (and hence by a
change of variable the entire two parameter Poisson Dirichlet class
discussed in Pitman~(1996)), classes of models based on the Pareto
and log-Normal distribution, and the class of Generalized Inverse
Gaussian~(GIG) models, among many others.

Let $N$ denote a Poisson random measure on the space
$(0,\infty)\times(0,\infty)\times\Xcr$, with mean intensity \Eq
\nu(dr,dv,dx)=\theta r^{-1}{\mbox e}^{-r}\Ucr(dv)H(dx)\label{GGC}
\EndEq where $\Ucr$ is a non-negative non-decreasing function on
satisfying, $\Ucr(0)=0$,\Eq \int_{0}^{1}|\ln(v)|\Ucr(dv)<\infty
{\mbox { and }}\int_{1}^{\infty}1/y\Ucr(dy)<\infty.
\label{Gcond}\EndEq Note importantly, that it is possible for
$\Ucr(\infty)=\infty$. In particular, the condition~\mref{Gcond}
is true if and only if
    \Eq
    \psi(\lambda)=\theta\,\int_{0}^{\infty}\int_{0}^{\infty}(1-{\mbox
    e}^{-\lambda(r/v)}) r^{-1}{\mbox
    e}^{-r}\Ucr(dv)=\theta\,\int_{0}^{\infty}\ln(1+\lambda/v)
    \Ucr(dv)<\infty.
    \label{GGC2}
    \EndEq
This allows one to define random variables
$T=G_{\theta}M_{\theta}=\int_{0}^{\infty}\int_{0}^{\infty}(r/v)N(dr,dv)$,
where $G_{\theta}$ is independent of $M_{\theta}$, and has L\'evy
exponent $\psi(\lambda)$. Let $\Tcr$ denote the family of all
random variables, $T$, such that the L\'evy exponent of $T$
satisfies~\mref{GGC2}. That is $\Tcr$ is the class of random
variables with distributions which are generalized gamma
convolutions. This may be extended by random variables $T+a$ for
$a\ge 0$. The quantity $\Ucr$ uniquely determines the distribution
of $T$, and we refer to $\Ucr$ as the {\it Thorin measure}.
Additionally,
    $$
    \mu(dx)=G_{\theta}\sum_{i=1}^{\infty}({\tilde Q}_{i}/V_{i})\delta_{Z_{i}}(dx)
    =\int_{0}^{\infty}\int_{0}^{\infty}(r/v)N(dr,dv,dx)
    $$
and $M_{\theta}=\sum_{i=1}^{\infty} ({\tilde Q}_{i}/V_{i})$. $\mu$
may be referred to as a GGC random measure. Then call
$P_{M_{\theta}}$ a GGC NRM if it has a representation as an
$h$-biased random probabilty measure, here $s:=(r,v)$, with
$h(r,v)=r/v$, given by
    \Eq
    P_{M_{\theta}}(dx)=\frac{\sum_{i=1}^{\infty} ({\tilde
    Q}_{i}/V_{i})\delta_{Z_{i}}(dx)}{M_{\theta}}=
    \frac{\int_{0}^{\infty}\int_{0}^{\infty}(r/v)N(dr,dv,dx)}{T}
    \label{meanNMR}
    \EndEq
where $(\tilde Q_{i})_i$ has the same marginal distribution as
in~\mref{DPrep}, but is now independent of the sequence $(Z_{i})_i$
of i.i.d. random variables whose distribution is $H$. Additionally,
both sequences are independent of $(V_{i})_i$. Note that the
distribution of the sequence is derived from the points of a Poisson
random measure with mean intensity $\Ucr$. In the special case where
$\Ucr$ is a probability measure (finite measure), the $(V_{i})_i$
are i.i.d. $\Ucr$. However this not always true. In fact the
obtainment of many interesting classes, such as the stable law,
require that $\Ucr$ is not a finite measure.

\Remark We mention that if $\Ucr$ is a finite measure then
$M_{\theta}\overset {d}=\int_{0}^{\infty}1/yD_{\theta}(dy)$,
$D_{\theta}$ being a Dirichlet process with shape parameter $\theta
\Ucr$. That is $M_{\theta}$ corresponds to a class of (positive)
Dirichlet mean functionals. That is, a subset of the Dirichlet
process functionals whose study was initiated by Regazzini and
Cifarelli~(1990). However $M_{\theta}$ constitutes a wider class of
positive random variables as the representation
$\sum_{i=1}^{\infty}{\tilde Q}_{i}\delta_{V_{i},Z_{i}}$ does not in
general correspond to a Dirichlet process unless the $V_{i}$'s are
iid. Rather $G_{\theta}\sum_{i=1}^{\infty}{\tilde
Q}_{i}\delta_{V_{i},Z_{i}}$ is a gamma process with possibly
sigma-finite shape measure $\theta \Ucr H$. As we shall see, this
generality of $M_{\theta}$ allows us much greater flexibility as
there are many cases where the distribution of $M_\theta$ and $T$
can be deduced. \EndRemark

\Remark Note interestingly, by independence of $G_{\theta}$ and
$M_{\theta}$, the fact that $G_{\theta}$ is gamma distributed
yields,
$$
E[{\mbox e}^{-uT}]=E[{\mbox
e}^{-uG_{\theta}M_{\theta}}]=\int_{0}^{\infty}{(1+uv)}^{-\theta}f_{M_{\theta}}(v)dv=
{\mbox e}^{-\theta\int_{0}^{\infty}\ln(1+u/y)\Ucr(dy)}
$$
where $f_{M_{\theta}}$ denotes the density of $M_{\theta}.$ This is
essentially the identity established in Cifarelli and
Regazzini~(1990) for Dirichlet mean functionals. Note that the
choice of $\Ucr$, uniquely determines the distribution of
$M_{\theta}$. Hence when $\Ucr(\infty)=\infty$, such distributions
are not captured by the current literature on Dirichlet mean
functionals. Note additionally that the theory of GGC has been
extended by Thorin~(1978) [see Bondesson~(1992)] to include
distributions on the entire real line. See also Lijoi and
Regazzini~(2004). Here of course we require that $T$ is positive.
\EndRemark

\subsection{Posterior Distribution of GGC NRM}
Now to establish the posterior distribution, first note that
$\mu_{n(\p)}|U_{n}=u,\X$ has mean intensity
$$
\nu_{u}(dr,dv,dx)=\theta r^{-1}{\mbox e}^{-r(1+u/v)}\Ucr(dv)H(dx)
$$
We recognize, from James~(2005a, Proposition 2.1),  that the
change in the intensity from $\nu$ to $\nu_{u}$ is due to
exponential tilting by ${\mbox e}^{-uT}$. It is useful to see
explicitly how this operation affects the Thorin measure, $\Ucr$,
and indeed how this affects the resulting distribution of
$M_{\theta}$.
\begin{prop} Let $T\in \Tcr$ defined by the Thorin measure
$\Ucr(dv):=\omega(v)dv$. Then suppose that $T_{b}$ is the random
variable with density $\propto {\mbox e}^{-bt}f_{T}(t)$ for some
$b>0$. Then it follows that $T_{b}\in \Tcr$, with Thorin measure
$\omega(z-b)I\{z>b\}dz$. This follows from the fact that its
L\'evy exponent is expressible as
$$
\psi_{b}(\lambda)=\int_{b}^{\infty}\ln(1+\lambda/z)\omega(z-b)dz=
\int_{0}^{\infty}\ln(1+\lambda/(y+b))\omega(y)dy.$$ Equivalently,
$T_{b}\overset {d}=G_{\theta}M_{\theta,b}=
\int_{0}^{\infty}\int_{0}^{\infty}(r/y)N_{b}(dr,dy)$, where
$N_{b}$ is a Poisson random measure with
$\E[N_{b}(dr,dz,dx)]=\theta r^{-1}{\mbox
e}^{-r}\omega(z-b)I\{z>b\}dzH(x)$.
$M_{\theta,b}:=\sum_{i=1}^{\infty}{\tilde Q}_{i}/(V_{i}+b)$ is the
mean functional whose law is induced by the tilting operation
\end{prop}
\Proof This follows from a straightforward change of variables or
using the fact that $\psi_{b}(\lambda)=\psi(\lambda+b)-\psi(b)$.
\EndProof

Now examining Proposition 7.1, with $b=U_{n}$, it follows that,
$$\mu_{n(\p)}\overset {d}=G_{\theta}\sum_{i=1}^{\infty}{\tilde
Q}_{i}{(V_{i}+U_n)}^{-1}\delta_{Z_{i}}.$$ Additionally
$T_{n(\p)}\overset {d}=G_{\theta}M_{\theta,U_{n}}$. Setting
$J_{j,n}:=(\Delta_{j,n},V_{j,n})$, its joint distribution given
$U_{n},\X$ is given by,
$$
\PR(\Delta_{j,n}\in dr,V_{j,n}\in dv|u)=\frac{{\mbox
e}^{-r(1+u/v)}r^{e_{j}-1}v^{-e_{j}}\theta
\Ucr(dv)dr}{\kappa_{e_{j}}(u)},
$$
where $
\kappa_{e_{j}}(u)=\theta\Gamma(e_{j,n})\int_{0}^{\infty}(v+u)^{-e_{j}}\Ucr(dv).
$ This shows that the distribution of each
$\Delta_{j,n}|V_{j,n},U_{n}$ is equivalent in distribution to
$G_{j,n}{(1+U_{n}/V_{j,n})}^{-1}$, where $G_{j,n}$ denotes a gamma
random variable with shape $e_{j}$ and scale $1$ independent of
$(U_{n},V_{j,n})$. The distribution of $V_{j,n}|U_{n},\X$ is given
by the density \Eq \PR(V_{j,n}\in
dv|u)=\frac{{(v+u)}^{-e_{j}}\theta\Gamma(e_{j})
\Ucr(dv)}{\kappa_{e_{j}}(u)}=\frac{{(v+u)}^{-e_{j}}
\Ucr(dv)}{\int_{0}^{\infty}{(y+u)}^{-e_{j}}\Ucr(dy)}
\label{denV}\EndEq Additionally, we use the fact that the prior
and posterior distribution of
$$U_{n}T=U_{n}G_{\theta}M_{\theta}=U_{n}\[G_{\theta}\sum_{i=1}^{\infty}{\tilde
Q}_{i}{(V_{i}+U_n)}^{-1}+\sum_{j=1}^{n(\p)}G_{j,n}{(V_{j,n}+U_{n})}^{-1}\]\overset
{d}=\Gamma_{n},$$ where $\Gamma_{n}$ is independent of $N$ and
$\X$. Furthermore notice that, \Eq
m_{n}(u)=E[T^{n}_{n(\p)}|u]=\frac{\Gamma(\theta+n)}{\Gamma(\theta)}
\frac{\int_{0}^{\infty}{(1+uv)}^{-(\theta+n)}v^{n}f_{M_{\theta}}(dv)}
{\int_{0}^{\infty}{(1+uv)}^{-\theta}f_{M_{\theta}}(dv)}\label{GGCmom}
\EndEq

These facts lead to a non-obvious description of the posterior
distribution given $U_{n},\X$.

\begin{thm}  Suppose that $P_{M_{\theta}}$, is the
NRM defined by ~\mref{meanNMR}. Then the following results hold
\Enumerate
\item[(i)]
The posterior distribution of $P_{M_{\theta}}|U_{n},\X$ is
equivalent to the distribution of the random probability measure
\Eq
U_{n}\frac{G_{\theta+n}}{\Gamma_{n}}\[\frac{G_{\theta}}{G_{\theta+n}}\sum_{i=1}^{\infty}{\tilde
Q}_{i}{(V_{i}+U_{n})}^{-1}\delta_{Z_{i}}+
\sum_{j=1}^{n(\p)}Q_{j,n}{(V_{j,n}+U_{n})}^{-1}\delta_{Y_{j}}\]
\label{postD2} \EndEq where, the gamma random variable
$G_{\theta+n}=G_{\theta}+\sum_{j=1}^{n(\p)}G_{j,n}$, and
independent of $U_{n},\Gamma_{n},(Z_{i})$, the $({\tilde Q}_{i})$
are equivalent in distribution to those in ~\mref{DPrep}, and are
independent of $\X$. Additionally, the vector
$(Q_{1,n},\ldots,Q_{n(\p),n})$ is independent of
$(U_{n},\Gamma_{n},(Z_{i}),G_{\theta},G_{\theta+n})$ and given
$\p$, is a Dirichlet $(e_{1},\ldots,e_{n(\p)};\theta)$ vector. The
relevant distributions of $V_{j,n},U_{n}$ are specified by
~\mref{denV} and~\mref{margUMp}.
\item[(ii)] Equivalently, one may write~\mref{postD2} as
    $$
    \frac{U_{n}T_{n(\p)}}{\Gamma_{n}}\sum_{i=1}^{\infty}{\tilde
Q}_{i,n}\delta_{Z_{i}}+ \[1-\frac{U_{n}T_{n(\p)}}{\Gamma_{n}}\]
\sum_{j=1}^{n(\p)}Q^{*}_{j,n}\delta_{Y_{j}}
    $$
where
$${\tilde Q}_{i,n}=\frac{{\tilde
Q}_{i}{(V_{i}+U_{n})}^{-1}}{\sum_{i=1}^{\infty}{\tilde
Q}_{l}{(V_{l}+U_{n})}^{-1}}\quad {\mbox { and }}\quad
Q^{*}_{j,n}=\frac{Q_{j,n}{(V_{j,n}+U_{n})}^{-1}}{\sum_{l=1}^{n(\p)}Q_{l,n}{(V_{l,n}+U_{n})}^{-1}}.$$
\EndEnumerate
\end{thm}
\subsection{$\Mcr$ and some connections to the Bondesson Class}
Before saying more about the distribution of $U_{n}, \X$, we next
describe an important subclass of GGC random variables which
interestingly are connected to the distribution of $U_{n}$. As in
Kent and Tyler~(2001, p.257) let $S$ denote a random variable on
$(0,\infty)$ with density
    \Eq
    f_{S}(s)=Cs^{\beta-1}\prod_{j=1}^{m}{(1+c_{j}s)}^{-\gamma_{j}}{\mbox
    { for }}s>0,
    \label{BondClass}
    \EndEq
where $C$ is a normalizing constant, $m\ge 1$, $\beta>0$,
$c_{j}>0$,$\gamma_{j}>0$.
Then the class containing the
densities ~\mref{BondClass}, together with their weak limits,
constitutes the Bondesson $\Bcr$ sub-class of GGC models. We write
$S\in \Bcr$ if the density of $S$ has the form
in~\mref{BondClass}. Note that in the non-limiting case, the
corresponding $\Ucr(\infty)=\beta<\infty$. That is, in this case,
the $(V_{i})$ are iid with distribution $\beta^{-1}\Ucr$. The
$\Bcr$ class contains the Stable distributions of index
$\alpha=1/k$ for $k=2,3\ldots$. The gamma distribution, Pareto,
Log-Normal random variables, generalized inverse gaussian, among
many others. [See Bondesson(1992, Chapter 5.6)]. This class is
known to be {\it hyperbolically completely monotone} and hence
{\it self-decomposable}. See Bondesson~(1992), Steutel and Van
Harn~(2004, Chapter 5) and Kent and Tyler~(2001, p.257) for
further details.
\begin{prop} Let
$U_{n}=\Gamma_{n}/T=\Gamma_{n}/(G_{\theta}M_{\theta})$, then its
distribution is described as follows. \Enumerate
\item[(i)] The density of $U_{n}|\p$ is, \Eq
f_{U_{n}}(u|\p)\propto
u^{n-1}\[\int_{0}^{\infty}{(1+uv)}^{-\theta}f_{M_{\theta}}(dv)\]
\prod_{j=1}^{n(\p)}\int_{0}^{\infty}(v+u)^{-e_{j}}\Ucr(dv)
\label{margUMp}\EndEq
\item[(ii)]The marginal density of $U_{n}$ is \Eq
\frac{u^{n-1}}{\Gamma(n)}\E[{\mbox
e}^{-uG_{\theta}M_{\theta}}{(G_{\theta}M_{\theta})}^{n}]
=u^{n-1}\frac{\Gamma(\theta+n)}{\Gamma(n)\Gamma(\theta)}\int_{0}^{\infty}{(1+uv)}^{-(\theta+n)}
v^{n}f_{M_{\theta}}(dv) \label{margUM}\EndEq This implies that
that for any integrable function $g$,
    $$
    \E[g(U_{n})]=\int_{0}^{1}\E\[g\(\frac{1-y}{yM_{\theta}}\)\]
    \Bcr(dy|\theta,n)
    $$
\EndEnumerate \end{prop}

\begin{prop} Suppose that $T\in \Tcr$, then
$U_{n}:=\Gamma_{n}/T=\Gamma_{n}/(G_{\theta}M_{\theta})$ has the
following properties. \Enumerate \item[(i)] If
$T:=G_{\theta}M_{\theta}$ $\in$ $\Bcr$, then $U_{n}\in \Bcr$
However, $T\in \Tcr$ does not imply that $U_{n}\in \Tcr.$
\item[(ii)] The distribution of $U_{n}|M_{\theta}=v$ is in $\Bcr$.
with density of the form in~\mref{BondClass}, with parameters,
$\beta=n$, $m=1$, $c_{1}=v$, $\gamma_{1}=\theta+n$. Equivalently,
$M_{\theta}U_{n}$ is a gamma-gamma density, and hence in $\Bcr$.
The form of the density coincides with $U_{n}|M_{\theta}=1.$ That
is $c_{1}=1$
\item[(iii)] The distribution of $U_{n}|M_{\theta}=v,(V_{j,n}=v_{j}),\p$ is in
$\Bcr$, with $\beta=n$, $m=n(\p)+1$,$c_{j}=1/v_{j}$,
$\gamma_{j}=e_{j}$, for $j=1,\ldots,n(\p)$ and $c_{n(\p)+1}=v$,
$\gamma_{n(\p)+1}=\theta$. Specifically the conditional density is
given by
$$
f_{U_{n}}(u|v,(v_{j}))=Cu^{n-1}{(1+uv)}^{-\theta}\prod_{j=1}^{n(\p)}{(1+u/v_{j})}^{-e_{j}}
$$
\item[(iv)] Let $Y_{n}=n/U_{n}$. Then $Y_{n}\in \Tcr$  and hence
this family of densities is dense in the class of all $\Tcr$.
\EndEnumerate
\end{prop}
\Proof Statement (i) is an immediate consequence of
Bondesson~(1992, Theorem 6.2.1, p.92.).  See also Kent and
Tyler~(2001, statement 7, p. 257). Statement (ii) and (iii) follow
from an augmentation and matching with~\mref{BondClass}.
Statement(iv) is read from Bondesson~(1992, p. 92).\EndProof

 The next result gives a form of the EPPF.
\begin{prop} Suppose that $\p$ denotes the partition derived from $P_{M_{\theta}}$. Then the following results hold.
\Enumerate
\item[(i)] The conditional distribution of
 $\p|U_{n}$ is given by \Eq PD(\p|\theta)\[
\frac{\int_{0}^{\infty}{(1+uv)}^{-\theta}f_{M_{\theta}}(dv)}
{\int_{0}^{\infty}{(1+uv)}^{-(\theta+n)}v^{n}f_{M_{\theta}}(dv)}\]
\prod_{j=1}^{n(\p)}\int_{0}^{\infty}(y+u)^{-e_{j}}\Ucr(dy)
\label{jointDDP2}\EndEq
\item[(ii)] The EPPF may be expressed as,
\begin{multline}
PD(\p|\theta)\frac{\Gamma(\theta+n)}{\Gamma(n)\Gamma(\theta)}\int_{0}^{\infty}u^{n-1}\[
{\int_{0}^{\infty}{(1+uv)}^{-\theta}f_{M_{\theta}}(dv)}\]\,\times\\
\[\prod_{j=1}^{n(\p)}\int_{0}^{\infty}(y+u)^{-e_{j}}\Ucr(dy)\]du
\label{jointDDP2}
\end{multline}
\EndEnumerate
\end{prop}

\Remark The result in~\mref{theile} combined with~\mref{GGCmom}
and the form of $\kappa_{j}$ leads to interesting relationships
between $\Ucr$ and $M_{\theta}$. In particular, using the property
that $\kappa_{1}(u)=m_{1}(u)$ leads to the following identity,
    $$
    \int_{0}^{\infty}(v+u)^{-1}\Ucr(dv) =
    \frac{\int_{0}^{\infty}{(1+uv)}^{-(\theta+1)}vf_{M_{\theta}}(dv)}
    {\int_{0}^{\infty}{(1+uv)}^{-\theta}f_{M_{\theta}}(dv)},
    $$
\EndRemark

\Remark The unified representation of $P_{M_{\theta}}$ and the
characterization of its posterior distribution given
in~\mref{postD2} in terms of the Dirichlet process $(Q_{i})$ has
many interesting implications. For instance, it suggests that one
can use a variant of the Blocked Gibbs algorithms in Ishwaran and
James~(2001, 2003b) to approximately sample realizations of
$P_{M_{\theta}}$ and its posterior process for many different
classes of models which are not Dirichlet processes. \EndRemark

\subsection{Some specific examples of GGC NRM}
In this section we now highlight some important specific cases.
First it is interesting to note the the study of the GGC is
primarily about establishing the fact that $T\in \Tcr$, and possibly
identifying $\Ucr$. In contrast, the study of Dirichlet process mean
functionals involves identification of the distribution of
$M_{\theta}$, when $\Ucr$ is a pre-specified finite measure. We see
these two approaches as complementary to one another. We point out
that explicit forms for the $\Ucr$ are not known in every case.
However, importantly there are many examples of $T$ which are known
to be in $\Tcr$. For instance, an explicit form of $\Ucr$ is not
known for the Pareto distribution. However, as we have shown, many
interesting applications involving sampling from $\Mcr$ and mixture
models can still be conducted if one knows $m_{n}(u)$ or the
cumulants $\kappa_{n}(u)$. Similarly, from the Dirichlet process
literature the explicit law of $M_{\theta}$ is not known in many
cases, however we can choose $\Ucr$ to be from a vast selection of
probability distributions. Thus as stated earlier, in that case, we
have an explicit description of the iid distribution of the sequence
$(V_{i})$.

As some general examples, one could choose $T=S\in \Bcr$, as
defined in \mref{BondClass}. Note also that $T=S^{q}\in \Tcr$ for
$|q|\ge 1$ [see Bondesson~(1979, Corollary 1)]. Here we give some
precise examples from the literature. \Remark Note carefully that
we only need to show that $T$ has a particular law to establish
the law of the NRM $P_{M_{\theta}}:=\mu/T$. This is due to the
fact that the Laplace functional of $\mu$ evaluated at some
bounded measurable functional $g$ is determined by
$$
-\log\E[{\mbox e}^{-\mu(g)}]=\int_{\mathcal X}\psi(g(x))H(dx)
$$
for $\psi(g(x))$ given by replacing $\lambda$ with $g(x)$
in~\mref{GGC2}\EndRemark

\subsubsection{Stable Case and related models}
As mentioned previously, the NRM based on the stable law have been
extensively studied by Pitman~(1996, 2002) and Pitman and
Yor~(1997). This class has numerous applications and also has a
tractable EPPF. The explicit posterior distribution of this class
was obtained by Pitman~(1996) by exploiting its explicit
stick-breaking representation. James (2002, section 5.3) gives an
alternative derivation working directly with the Levy measure of a
stable law. Here, we show that $P_{M_{\theta}}$ offers another
representation of the stable law NRM and hence alternative approach
to its analysis. We give some details of its posterior analysis
which are inherent to its representation in terms of a GGC NRM.
Further details can be deduced easily from the specific analysis of
these models in James~(2002, sections 5.3 and 5.4).

From Bondesson~(1992, p. 35), it is now easy to see that $T$ is
stable, and hence $P_{M_{\theta}}$ is a stable NRM,  if
$$
\Ucr_{\alpha}(dv)=\frac{\alpha}{\Gamma(1-\alpha)\Gamma(\alpha)}v^{\alpha-1}dv
{\mbox { for }} v>0.
$$
As checks, one can make the change of variable $z=r/v$,
in~\mref{GGC2} and integrate with respect to $\Ucr$ first. That is
for all $\theta>0$,
$\psi(\lambda)=C_{\alpha,\theta}\lambda^{\alpha}$ for some
constant $C_{\alpha,\theta}$. Note interestingly that, conditional
on $U_{n}$, we have
    \Eq
    \left(\frac{V_{j,n}}{U_{n}}+1\right)^{-1}\sim\Bcr(\e_{j}-\alpha, \alpha)
    \label{StableV}
    \EndEq
and, hence,
    $$
    G_{j,n}\left(\frac{V_{j,n}}{U_{n}}+1\right)^{-1}
    \sim\Gcr(e_{j}-\alpha)
    $$
for $j=1,\ldots, n(\p)$. That is, these quantities are independent
of $U_{n}$. Additionally, $\kappa_{e_{j}}(u)=\theta\alpha
u^{\alpha-e_{j}}{\Gamma(e_{j}-\alpha)}/\Gamma(1-\alpha)$. Ignoring
scale parameters it follows that the distribution of
$L_{n}:=U^{\alpha}_{n}|\p$ is $\Gcr(n(\p))$. One then establishes
that $\mu_{n(\p)}|U_{n},\p$ is a generalized gamma process, whose
Thorin measure is given by ~\mref{gengam} below, with $b=U_{n}$.
Moreover the distribution of $L^{1/\alpha}_{n}\mu_{n(\p)}|L_{n},\p$
is determined by the L\'evy measure
    \Eq
    \frac{L_{n}\alpha}{\Gamma(1-\alpha)\Gamma(\alpha)}r^{-1}{\mbox
    e}^{-r}{(v-1)}^{\alpha-1}I\{v>1\}dr dv H(dx).
    \label{Lmix}
    \EndEq
In particular this implies that the distribution of
$L^{1/\alpha}_{n}T_{n(\p)}|\p$ is $\Gcr(n(\p)\alpha)$. These formula
can be used in Proposition 7.1 and 7.3, to establish the known
results about the posterior distribution. Further details can be
deduced from James~(2002, section 5). The two parameter
Poisson-Dirichlet distribution with parameters $(\alpha,q)$ for
$0<\alpha<1$ and $q>-\alpha$ arises from the weighted law $\propto
T^{-q}\Pe(dN|\nu)$ as described in Proposition 8.2. It remains to
note that $L_{n,q}:=U^{\alpha}_{n,q}$ given $\p$ is
$\Gcr(n(\p)+q\alpha)$

\Remark Note also that, in the case where $T$ is the stable law,
the generalized Cauchy-Stieltjes transform of $M_{\theta}$ is
$$
\int_{0}^{\infty}{(1+\lambda
v)}^{-\theta}f_{M_{\theta}}(dv)={\mbox
e}^{-C_{\alpha,\theta}\lambda^{\alpha}}
$$
which can be easily inverted. This, in some sense, easiest example
is outside of the scope of the current theory of Dirichlet process
mean functionals as $\Ucr(\infty)=\infty$. \EndRemark

\subsubsection{Generalized Gamma} The class of generalized gamma process defined for $0<\alpha<1$,
and $b>0$ [see Brix~(1999)] arises from the tilting by ${\mbox
e}^{-bS_{\alpha}}$, where $S_{\alpha}$ is stable law of index
$\alpha$. The simplest case is when $\alpha=1/2$, where the
corresponding $T$ has an Inverse Gaussian distribution. Proposition
7.1 shows that the Thorin measure is given, in this case, by
    \Eq
    \frac{\alpha}{\Gamma(1-\alpha)\Gamma(\alpha)}{(v-b)}^{\alpha-1}dv
    \qquad\quad {\mbox { for }} \quad v>b.
    \label{gengam}
    \EndEq
The Thorin measure of $\mu_{n(\p)}|U_{n},\p$ is of the same form
as~\mref{gengam} with $b$ replaced by $U_{n}+b$, and hence is a
generalized gamma process. It follows that setting
$L_{n}=C_{\alpha,\theta}(U_{n}+b)^{\alpha}$, the distribution of
$L^{1/\alpha}_{n}\mu_{n(\p)}|L_{n},\p$ is the same as that for the
stable law determined by~\mref{Lmix}. Similar to~\mref{StableV} one
has, conditional on $U_{n}$,
    $$
    \left(\frac{V_{j,n}-b}{U_{n}+b}+1\right)^{-1}\sim \Bcr(\e_{j}-\alpha,
    \alpha)\quad {\mbox { and }}\quad
    G_{j,n}\left(\frac{V_{j,n}-b}{U_{n}+b}+1\right)^{-1}\sim\Gcr(e_{j}-\alpha)
    $$
for $j=1,\ldots, n(\p)$. With $\kappa_{e_{j}}(u)=\theta\alpha
{(u+b)}^{\alpha-e_{j}}{\Gamma(e_{j}-\alpha)}/\Gamma(1-\alpha)$. The
density of $U_{n}|\p$ is
$$
f_{U_{n}}(u|\p)\propto {(u+b)}^{n(\p)\alpha-n}u^{n-1}{\mbox
e}^{-C_{\alpha,\theta}[{(u+b)}^{\alpha}-b^{\alpha}]}
$$
Using a Binomial expansion, the distribution of $L_{n}|\p$ is
given, for all $b>0$, by
    $$
    f_{L_{n}}(w|\p)=\frac{\sum_{k=0}^{n-1}{n-1\choose
    k}{(-1)}^{k}{w}^{-k/\alpha}I\{w>C_{\alpha,\theta}\}\Gcr(w|n(\p))}{\sum_{k=0}^{n-1}{n-1\choose
    k}{(-1)}^{k}\E[G^{-k/\alpha}_{n(\p)}I\{G_{n(\p)}>C_{\alpha,\theta}\}]
    }
    $$
The normalizing constant can be used to yield an explicit
expression for the EPPF.

\Remark Note that the for the range $\alpha=0$ and $b>0$, the
generalized gamma process equates with the the gamma process.
\EndRemark

\subsubsection{Generalized Inverse Gaussian}
A more challenging class is the Generalized Inverse Gaussian (GIG)
class of models. First set $\theta=1$. Let $\lambda$, $v$ and
$\delta$ be such that $0<\lambda<\infty$, while $v$ and $\delta$
are non-negative and not simultaneously $0$. As in
Barndorff-Nielsen and Shephard (2001), $T$ is GIG$(\lambda,
\delta, v)$ if its density is of the form
  \Eq
  f_{T}(t|\lambda,\delta,v)=\frac{{({v/\delta})}^{\lambda}}{2K_{\lambda}(\delta
  v)}t^{\lambda-1}\exp\{-\frac{1}{2}(\delta^{2}t^{-1}+v^{2}t)\}
  \EndEq
where $K_{\lambda}$ is a Bessel function of the third kind. When
$\delta=0$ and $\lambda>0$, $v>0$, GIG$(\lambda,0,v)$ equates with
the gamma distribution.  When $\lambda<0$, $\delta>0$ and $v=0$,
then GIG$(\lambda, \delta,0)$ is a reciprocal, or inverse gamma
distribution. Using the parametrization, $\lambda=-a$, for $a>0$,
and $b=\delta^{2}/2$, yields the density of an inverse gamma
distribution with parameters, $a,b$,with density
    \Eq
    \frac{{(2/\delta^{2})}^{\lambda}}{\Gamma(-\lambda)}
    t^{\lambda-1}\exp\{-\frac{1}{2}(\delta^{2}t^{-1})\}
    =\frac{{b}^{a}}{\Gamma(a)}t^{-a-1}\exp\{-b t^{-1}\}
    \label{recip}
    \EndEq
A special case of this is when $\lambda=-1/2$ leading to a stable
law of index $1/2$. The Inverse Gaussian distribution defined by
setting $\lambda=-1/2$, $\delta>0$, and $v>0$ that is a
GIG$(-1/2,\delta,v)$. The hyperbolic distribution coincides with the
case of $\lambda=1$. Now define,
    \Eq
    g_{|\lambda|}(x)=\frac{2}{\pi^{2}}\frac{1}{x[J^{2}_{|\lambda|}
    (\sqrt{x})+N^{2}_{|\lambda|}(\sqrt{x})]}
    \label{www}
    \EndEq
where $J_{v}$ and $N_{v}$ are Bessel functions of the first and
second kind respectively. The expression~\mref{www} is central to
a body of work on the infinite divisibility of student
t-distribution and generalized gamma convolutions. One has for
$m=0,1,2\ldots,$
    \Eq
    g_{m+1/2}(x)=\frac{2}{\pi^{2}}\frac{x^{(m-1)/2}}{\prod_{i=1}^{m}(x+a^{2}_{j})}
    \label{integercase}
    \EndEq
where the $a_{1},\ldots,a_{m}$ are the zeros of $K_{m+1/2}(z)$. It
is known that $\Ucr$ is given by
    \Eq
    \Ucr(\ddr x)=I_{\{x\ge v^{2}/2\}}
    \[\delta^{2}\int_{v^{2}/2}^{x}g_{|\lambda|}(2{\delta}^{2}y-\delta^{2}v^{2})dy
    +\max(0,\lambda)\]\,\ddr x
    \EndEq
The simplest cases correspond to the gamma distribution, that is the
Dirichlet process, and cases covered by ~\mref{integercase}. Setting
$m=0$, $v=0$ in~\mref{integercase} coincides with the $\Ucr_{1/2}$
of a stable $(1/2)$ law. When $v>0$, one obtains the
Inverse-Gaussian distribution. Setting $\lambda=-3/2$ gives $m=1$,
now with $v>0$,$\delta>0$ corresponds yields to the Thorin measure
given by
    $$
    \Ucr(\ddr y)=\frac{2}{\pi^{2}}{(2y-v^{2}+a^{2}_{1}\delta^{-2})}^{-1}I\{y>v^{2}/2\}\, \ddr y
    $$
The simplest case arises if one further sets
$v^{2}=a^{2}_{1}\delta^{-2}$. The other cases
involving~\mref{integercase} are a slightly more complex but
certainly can be handled. For the general case, using~\mref{www}
one can calculate the $\kappa_{n}$ from the moments $m_{n}$ which
are obtained as ratios of Bessel functions $K_{\lambda}$. Making
the substitution $u=(w^{2}-v^{2})/2$ for $w\ge v>0$, gives the
Laplace transform and the $m_{n}$ as follows,
$$
\phi((w^{2}-v^{2})/2)=\frac{v^{\lambda}K_{\lambda}(\delta
w)}{w^{\lambda}K_{\lambda}(\delta v)}\quad  {\mbox { and }}\quad
m_{n}((w^{2}-v^{2})/2)=\delta^{n}w^{-n}\frac{K_{n+\lambda}(\delta
w)}{K_{\lambda}(\delta w)}
$$
The marginal density of $L_{n}=\sqrt{(2U_{n}+v^{2})}$ is given by,
$$
f_{L_{n}}(w)=v^{\lambda}{2}^{-(n-1)}{(w^{2}-v^{2})}^{n-1}w^{-(n+\lambda-1)}
\delta^{n}\frac{K_{n+\lambda}(\delta
w)}{K_{\lambda}(\delta v)}I\{w\ge v\}
$$
Note that one can use the further simplification for
$n=0,1,2\ldots$
$$
K_{n+1/2}(z)=\sqrt{\frac{\pi}{2x}}{\mbox
e}^{-x}\sum_{k=0}^{n}\frac{(n+k)!}{2^{k}(n-k)!k!}x^{-k}
$$
Details may be deduced from Barndorff-Nielsen and Shephard (2001).
\Remark We mention again, that although the exact form of the
density for $L_{n}$ or $U_{n}$ appears to be complicated, these are
easily simulated using the fact that $U_{n}=\Gamma_{n}/T$. For
example in the case of the inverse gamma distribution
$U_{n}=\Gamma_{n}G_{a}$. The exact representation of the density is
of interest for instance in possible connections and interpretations
to the theory of special functions. Lijoi and Regazzini~(2004) is an
example of recent work exploring the interface between special
functions and problems arising in Bayesian nonparametrics. See also
James (2005b). \EndRemark
\subsubsection{First passage time distribution}
The next example, taken from Bondesson~(1992, p.37), involves
$\Ucr$ which is a proper distribution. For simplicity set
$\theta=1$. Let $1/2\leq p<1$, then $T$ has a {\it first passage
time} distribution if its moment generating function evaluated at
$\lambda$, has the form for $b=2\sqrt{p(1-p)}\leq 1$,
$$
\frac{1-\lambda-\sqrt{(1-\lambda)^{2}-b^{2}}}{2(1-p)}
$$
In this case, $
\Ucr(dy)=\frac{1}{\pi}{(y-1+b)}^{-1/2}{(1+b-y)}^{-1/2}dy,{\mbox {
for }}1-b<y<1+b.$
\subsubsection{Some examples from the Dirichlet process mean
functional} As mentioned previously, when  $\Ucr$ is finite then
$M_{\theta}$ has the law of a Dirichlet process mean functional.
That is, taking $\Ucr$ as a distribution function, the $(V_{i})$ are
iid $\Ucr$. Due to the work of Cifarelli and Regazzini~(1990), the
law of $M_{\theta}$ is known to often have a complex density which
is not commonly seen in the literature. It is of course a simple
matter to then obtain an expression for the distribution of
$T=G_{\theta}M_{\theta}.$ Here we state two examples. First, suppose
that $\theta=1$ and $1/V_{i}$ is chosen to be a uniform distribution
on $(0,1)$. Then it is known that the distribution of
$M_{1}=\int_{0}^{1}(1/y)D_{1}(dy)$ has a density given by
$$
f_{M_{1}}(v)=\frac{\mbox e}{\pi}v^{-v}{(1-v)}^{-(1-v)}\sin(\pi
v){\mbox { for }}0<v<1
$$
That is $\Ucr(dy)=y^{-2}dy$ for $y>1$. The final example may be
found in Cifarelli and Mellili~(2000). Suppose that $1/V_{i}$ is
$\Bcr(1/2,1/2)$, that is the Arc-sine law. Then for all
$\theta>0$, $M_{\theta}$ is $\Bcr(\theta+1/2,\theta+1/2)$.

\section{Appendix}
\subsection{Proof of Theorem 2.1, Propositions 2.1-2.3}
\Proof An intial description of the posterior distribution of
$N|\X$, follows as a simple variant of Theorem 3.2 in
James(2005a). First note that the result in James~(2005a, Theorem
3.2) holds obviously with $h(s)$ in place of $s$. One can easily
verify this by using James~(2005a, Theorem 3.1). Then by using
that result with $h(s,N):=h(s)/T$, it follows that the posterior
distribution of $N|\X$ is equivalent to the distribution of the
random measure $N^{*}_{n}={\tilde
N}+\sum_{j=1}^{n(\p)}\delta_{J_{j,n},Y_{j}}$, where the joint law
of ${\tilde N},(J_{j,n})|\X$ is proportional to the joint measure,
\Eq \frac{1}{{({\tilde
T}+\sum_{j=1}^{n(\p)}h(J_{j,n}))}^{n}}\Pe(d{\tilde
N}|\nu)\prod_{j=1}^{n(\p)}{[h(J_{j,n
})]}^{e_{j}}\rho(dJ_{j,n}|Y_{j}) \label{basicpost}\EndEq Note that
${\tilde N}$, ${\tilde T}$, corresponds to $N_{n(\p)}, T_{n(\p)}$
and $T={\tilde T} +\sum_{j=1}^{n(\p)}h(J_{j,n})$. Additionally
$\Pe(d\tilde N|\nu)$ is a Poisson law with intensity $\nu$, which
importantly is the same as the prior law of $N$. That is, under
$\Pe(d\tilde N|\nu)$,  $\Pe({\tilde T}\in dt)=f_{T}(t)dt$. Now an
application of the gamma identity~\mref{gammaid} yields a
posterior distribution of ${\tilde N},(J_{j,n}),U_{n}|\X$
proportional to, \Eq u^{n-1} {\mbox e}^{-u{\tilde T}}\Pe(d{\tilde
N}|\nu)\prod_{j=1}^{n(\p)}{[h(J_{j,n })]}^{e_{j}}{\mbox
e}^{-uh(J_{j,n})}\rho(dJ_{j,n}|Y_{j}) \label{postlast}\EndEq The
result then follows by applications of Bayes rule to
~\mref{postlast}. In particular notice that from Proposition 2.1.
of James~(2005a),
$$
{\mbox e}^{-uT}{\mbox e}^{\psi(u)}\Pe(dN|\nu)=\Pe(dN|\nu_{u})
$$
The description of $\Mcr$ appearing in Proposition 2.3 is an
immediate consequence of Theorem 3.2 of James~(2005a) combined
with the gamma identity.\EndProof

\subsection{The prediction rule} From Propositions 2.2 and 2.3 one can derive the
Bayesian prediction rule, i.e.
    \begin{align*}
    \Pe(X_{n+1}\in dx_{n+1}|\X)
    &=
    \frac{\Mcr(dX_{1},\ldots,d{x_{n+1}})}{\Mcr(dX_{1},\ldots,d{X_{n}})}\\
    &=
    \int_{0}^{\infty}\E\left[P^{*}_{n}(dx_{n+1})|u,\X\right]\,
    f_{U_n}(u|\X)\,\ddr u
    \end{align*}
One can rewrite the predictive distribution as a linear
combination of the measure $\eta$ and of a weighted empirical
distribution as follows
   \begin{equation}
    \label{predict}
    \Pe(X_{n+1}\in dx_{n+1} | \X)=
    \zeta^{(n)}\,\eta(dx_{n+1})
    +\frac{1}{n}\,\sum_{j=1}^{n(\p)}
    \zeta_j^{(n)}\,\delta_{Y_j}(dx_{n+1})
    \end{equation}
where
    $
    \zeta^{(n)}=\frac{1}{n}\,
    \int_0^{+\infty}\tau_{1}(u|x_{n+1}) \, f_{U_n}(u|\X)\,\ddr u
    $
and, for each $j=1, \ld, n(\p)$,
    $$
    \zeta_j^{(n)}=\int_0^{+\infty}\frac{\tau_{e_{j,n}+1}(u|Y_j)}{\tau_{e_{j,n}}(u|Y_j)}\,
    f_{U_n}(u|\X)\,du=\int_0^{+\infty}\E\[h(J_{j,n})|u,\X\]\,
    f_{U_n}(u|\X)\,du.
$$
See also James (2002) and Pr\"unster (2002). One immediately
notices that, in general, the empirical distribution
$\sum_{j=1}^{n(\p)}e_{j,n}\delta_{Y_j}/n$ is no longer sufficient
for prediction. This is in contrast to what happens with the
Dirichlet process.

\subsection{Results for weighted Poisson laws} Suppose that $g$ is
a positive measurable function, such that, without loss of
generality, $\E[g(N)]=\int_{\M}g(N)\Pe(dN|\nu)=1$. In this section
we describe what happens when $P$ is governed by a weighted
Poisson law $\Pe_{g}(dN|\nu)\propto g(N)\Pe(dN|\nu)$. We also
highlight the case where $g(N)=T^{-q}/\E[T^{-q}]$ for some
$-\infty<q<\infty.$
\begin{prop} Suppose that $P$ is governed by the weighted Poisson
law $\Pe_{g}(dN|\nu)$ described above. Then it follows that the
law of $(N_{n(\p)},(J_{j,n}),\X,U_{n})$ is proportional to
$$g(\tilde
{N}+\sum_{j=1}^{n(\p)}\delta_{s_{j},Y_{j}})\Pe(d\tilde
{N}|\nu_{u})\[\prod_{j=1}^{n(\p)}\Pe(J_{j,n}\in
ds_{j}|u)\]\Mcr(d\X|u)f_{U_{n}}(u)$$ where the $N_{n(\p)}:=\tilde
{N}$, and otherwise the laws above correspond to those given in
Theorems 2.1 and 3.1. An application of Bayes rule yields the
relevant marginal and posterior distributions.\end{prop} We now
describe an important special case of Proposition 8.1.

\begin{prop} Suppose that $w_{q}=\E[T^{-q}]<\infty$ for some
$-\infty<q<\infty$. Furthermore assume that $n+q\ge 0$, then the
NRM $P$ defined by the weigthed Poisson measure $\propto
T^{-q}P(dN|\nu)$, has the following properties. \Enumerate
\item[(i)]
The law of $(N_{n(\p)},(J_{j,n}),\X,U_{n}=u)$ is
    $$
    \Pe(d\tilde
    {N}|\nu_{u})\[\prod_{j=1}^{n(\p)}\Pe(J_{j,n}\in
    ds_{j}|u)\]\Mcr(d\X|u)f_{U_{n,q}}(u)
    $$
where
$f_{U_{n,q}}(u)=u^{q}[\Gamma(n)/\Gamma(n+q)]w^{-1}_{q}f_{U_{n}}(u)$.
\item[(ii)] The  density $f_{U_{n,q}}$ corresponds to a random variable
$U_{n,q}\overset {d}=\Gamma_{n+q}/T$, where $\Gamma_{n+q}$ is a
gamma random variable independent of $N$, and $T$ now has marginal
density $t^{-q}w^{-1}_{q}f_{T}(t)$.
\item[(iii)] The posterior distributions of $N$, and hence
$P$ and $\mu$, given $U_{n,q},\X$ is the same as in Theorem 2.1.
The marginal distributions of $\X$ and $\p$ are given by
$\int_{0}^{\infty}\Mcr(d\X|u)f_{U_{n,q}}(u)du.$
\EndEnumerate
\end{prop}
\Proof The result follows by applying the gamma identity to
$T^{-(n+q)}$.

\bigskip
\bigskip
{\Heading Acknowledgements} We would like to thank Lennart
Bondesson for a helpful conversation and references related to the
{\it real inversion formula}.

\begin{center}
\textsc{References}
\end{center}

\begin{footnotesize}
\noindent \textsc{Aldous, D, {\rm and} Pitman, J.} (2002).
\textit{Two recursive decompositions of Brownian bridge.} To
appear in S\'eminaire de Probabilit\'es XXXIX (2002). Available at
\texttt{arxiv.math.PR/0402399}.

\noindent \textsc{Arratia, R.A. , Barbour, A.D. {\rm and} Tavar\'e,
S.} (2003). \textit{Logarithmic combinatorial structures: a
probabilistic approach.} EMS Monographs in Mathematics. European
Mathematical Society (EMS), Zurich.

\noindent \textsc{Arratia, R.A., Barbour, A. D., {\rm and} Tavar\'e,
S.} (1999). The Poisson-Dirichlet distribution and the
scale-invariant Poisson process  \emph{Combin. Probab. Comput.}
\textbf{8}, 407-416.

\noindent \textsc{Barndorff-Nielsen, O.E.} and \textsc{Shephard, N.}
(2001). Non-Gaussian Ornstein-Uhlenbeck-based models and some of
their uses in financial economics. \textit{J. R. Stat. Soc. Ser. B
Stat. Methodol.}, {\bf 63}, 167-241.

\noindent \textsc{Bondesson, L.} (1979). A general result on
infinite divisibility \emph{Ann. Probab.} \textbf{7} 965-979.

\noindent \textsc{Bondesson, L.} (1992). Generalized gamma
convolutions and related classes of distributions and densities.
Lecture Notes in Statistics, 76. Springer-Verlag, New York.

\noindent \textsc{Cifarelli, D.M., {\rm and} Melilli, E.} (2000).
Some new results for Dirichlet priors. \emph{Ann. Statist.}
\textbf{28}, 1390--1413.

\noindent \textsc{Cifarelli, D. M., {\rm and} Regazzini, E.} (1990).
Distribution functions of means of a Dirichlet process. \emph{Ann.
Statist.} \textbf{18} (1990), 429--442.

\noindent \textsc{Derrida, B.} (1981). Random-energy models: an
exactly solvable model of disordered systems. \textit{Phys. Rev.
B}, \textbf{24}, 2613--2626.

\noindent  \textsc{Doksum, K.} (1974). Tailfree and neutral random
probabilities and their posterior distributions. {\it Ann.
Probab.}, {\bf 2}, 183--201.

\noindent  \textsc{Donnelly, P. {\rm and} Grimmett, G.} (1993). On
the asymptotic distribution of large prime factors. \textit{J.
London Math. Soc.}, \textbf{47}, 395--404.

\noindent  \textsc{Engen, S.} (1978). \textit{Stochastic abundance
models.} Chapman \& Hall, London.


\noindent \textsc{Escobar, M.D.} (1994). Estimating normal means
with the Dirichlet process prior \textit{J. Amer. Statist. Assoc.}
\textbf{89 } 268-277.

\noindent \textsc{Escobar, M.D. {\rm and} West, M.} (1995).
Bayesian density estimation and inference using mixtures.
\textit{J. Amer. Statist. Assoc.} \textbf{90}, 577--588.

\noindent \textsc{Ewens, W. J.} (1972). The sampling theory of
selectively neutral alleles. \emph{Theor. Popul. Biol.} \textbf{3}
87-112.

\noindent  \textsc{Ewens, W.J. {\rm and} Tavar\'e, S.} (1995). The
Ewens sampling formula. In \textit{Multivariate discrete
distributions} (Johnson, N.S., Kotz, S. and Balakrishnan, N.,
eds.). Wiley, New York.

\noindent \textsc{Feller, W.} (1971). \textit{An introduction to
probability theory and its applications.} Vol. II. Second edition
John Wiley and Sons, Inc., New York-London-Sydney

\noindent \textsc{Ferguson, T.S.} (1973). A Bayesian analysis of
some nonparametric problems. \textit{Ann. Statist.} \textbf{1},
209--230.

\noindent  \textsc{Freedman, D.A.} (1963). On the asymptotic
behavior of Bayes' estimates in the discrete case. \textit{Ann.
Math. Statist.}, \textbf{34}, 1386--1403.

\noindent  \textsc{Grote, M. {\rm and} Speed, T.P.} (2002).
Approximate Ewens formulae for symmetric over-dominance selection.
\textit{Ann. Appl. Prob.}, \textbf{12}, 637--663.

\noindent  \textsc{Hjort, N.L.} (1990). Nonparametric Bayes
estimators based on beta processes in models for life history
data. \textit{Ann. Statist.}, \textbf{18}, 1259--1294.

\noindent \textsc{Hjort, N.L.} (2003). Topics in non-parametric
Bayesian statistics. In Green, P.J. et al., editors, {\em Highly
Structured Stochastic Systems}. Oxford Statistical Science Series,
\textbf{27}, 455--478. Oxford University Press, Oxford.

\noindent  \textsc{Ishwaran, H.} and \textsc{James, L. F.} (2001).
Gibbs sampling methods for stick-breaking priors. \textit{J. Amer.
Stat. Assoc.}, \textbf{96}, 161--173.

\noindent  \textsc{Ishwaran, H.} and \textsc{James, L.F.} (2003a).
Generalized weighted Chinese restaurant processes for species
sampling mixture models. \textit{Statist. Sinica}, {\bf 13},
1211--1235.

\noindent \textsc{Ishwaran, H. {\rm and} James, L. F.}  (2003b).
Some further developments for stick-breaking priors: finite and
infinite clustering and classification \emph{Sankhy\=a} \textbf{65}
577-592.

\noindent \textsc{James, L.F.} (2002). Poisson Process Partition
Calculus with applications to exchangeable models and Bayesian
Nonparametrics. Manuscript. \\ Available at
\texttt{arxiv.math.PR/0205093}.

\noindent \textsc{James, L.F.} (2005a). Poisson process partition
calculus with applications to Bayesian L\'evy moving averages and
shot--noise processes.  \emph{Ann. Statist.}, to appear.\\
Available at \texttt{http://ihome.ust.hk/$\sim$lancelot/}

\noindent \textsc{James, L.F.} (2005b). Functionals of  Dirichlet
processes, the Cifarelli-Regazzini identity and Beta-Gamma
processes. \emph{Ann. Statist.}, to appear.

\noindent \textsc{James, L.F.} (2005c). Poisson calculus for spatial
neutral to the right processes. {\it Ann. Statist.}, to appear.

\noindent\textsc{Janson, S.} (2001). Asymptotic distribution for
the cost of linear probing hashing. \textit{Random Structures
Algortihms}, \textbf{19}, 438--471.

\noindent \textsc{Kent, J.T. {\rm and} Tyler, D.E.}(2001).
Regularity and uniqueness for constrained $M$-estimates and
redescending $M$-estimates. \emph{Ann. Statist.} \textbf{29}
252--265.


\noindent  {\sc Kingman, J.F.C.} (1975). Random discrete
distributions. {\it J. Roy. Statist. Soc., Series B}, {\bf 37},
1--22.

\noindent \textsc{Kim, Y.} (1999). Nonparametric Bayesian estimators
for counting processes. \emph{Ann. Statist.} \textbf{27} 562-588.

\noindent \textsc{Kolchin, V. F.}(1986) \textit{Random mappings.}
Translated from the Russian. With a foreword by S. R. S. Varadhan.
Translation Series in Mathematics and Engineering. Optimization
Software, Inc., Publications Division, New York.

\noindent  \textsc{Lavine, M.} (1992). Some aspects of P\'olya
tree distributions for statistical modelling. \textit{Ann.
Statist.} \textbf{20}, 1222--1235.

\noindent \textsc{Lijoi, A., Mena, R.H. {\rm and} Pr\"unster, I.}
(2005). Hierarchical mixture modelling with normalized inverse
Gaussian priors. \textit{J. Amer. Stat. Assoc.}, to appear.

\noindent \textsc{Lijoi, A. {\rm and} Regazzini, E.} (2004). Means
of a Dirichlet process and multiple hypergeometric functions.
\emph{Ann. Probab.} \textbf{32} 1469--1495.

\noindent \textsc{Liu, J.S.} (1996). Nonparametric hierarchichal
Bayes via sequential imputations. \textit{Ann. Statist.}
\textbf{24}, 911--930.

\noindent \textsc{Lo, A.Y.} (1984). On a class of Bayesian
nonparametric estimates: {I} Density estimates. \textit{Ann.
Statist.} \textbf{12}, 351--357.

\noindent \textsc{Lo, A.Y., Brunner, L.J. and Chan, A.T.} (1996).
Weighted Chinese restaurant processes and Bayesian mixture model.
Research Report Hong Kong University of Science and
Technology.\\Available at
\texttt{http://www.erin.utoronto.ca/~jbrunner/papers/wcr96.pdf.}

\noindent \textsc{Mauldin, R.D., Sudderth, W.D. {\rm and}
Williams, S.C.} (1992). P\'olya trees and random distributions.
\textit{Ann. Statist.} \textbf{20}, 1203--1221.

\noindent  \textsc{McCloskey, J.W.} (1965). A model for the
distribution of individuals by species in an environment.
\textit{Ph.D. thesis}, Michigan State University.

\noindent \textsc{Nieto-Barajas, L.E., Pr\"unster, I.} {\rm and}
\textsc{Walker, S.G.} (2004). Normalized random measures driven by
increasing additive processes. {\it Ann. Statist.} \textbf{32},
2343--2360.

\noindent  \textsc{Perman, M., Pitman, J., {\rm and} Yor, M.}
(1992). Size-biased sampling of Poisson point processes and
excursions. \textit{Probab. Theory Related Fields} \textbf{92},
21--39.

\noindent \textsc{Pitman, J.} (1995). Exchangeable and partially
exchangeable random partitions. \textit{Probab. Theory Related
Fields} \textbf{102}, 145--158.

\noindent  \textsc{Pitman, J.} (1996). Some developments of the
Blackwell-MacQueen urn scheme. In \textit{Statistics, Probability
and Game Theory. Papers in honor of David Blackwell} (Eds.
Ferguson, T.S., \textit{et al.}). Lecture Notes, Monograph Series,
\textbf{30}, 245--267. Institute of Mathematical Statistics,
Hayward.

\noindent  \textsc{Pitman, J.} (1997). Partition structures derived
from Brownian motion and stable subordinators. \textit{Bernoulli}
\textbf{3}, 79--96.

\noindent  \textsc{Pitman, J.} (2002). \textit{Combinatorial
stochastic processes.} Lecture notes for St.~Flour Summer School.

\noindent \textsc{Pitman, J.} (2003). Poisson-Kingman partitions.
In Goldstein, D.R., editor, {\it Science and Statistics: A
Festschrift for Terry Speed}. Lecture Notes, Monograph Series,
\textbf{40}, 1--35. Institute of Mathematical Statistics, Hayward.

\noindent  \textsc{Pitman, J.} and \textsc{Yor, M.} (1997). The
two-parameter Poisson-Dirichlet distribution derived from a stable
subordinator. \textit{Ann. Prob.} \textbf{25}, 855--900.

\noindent \textsc{Pr\"unster, I.} (2002). {\it Random probability
measures derived from increasing additive processes and their
application to Bayesian statistics}. Ph.d dissertation, University
of Pavia.

\noindent  \textsc{Regazzini, E., Lijoi, A. {\rm and} Pr\"unster,
I.} (2003).  Distributional results for means of random measures
with independent increments. \textit{Ann. Statist.} \textbf{31},
560--585.

\noindent  \textsc{Ruelle, D.} (1987). A mathematical
reformulation of Derrida's REM and GREM. \textit{Comm. Math.
Phys.}, \textbf{108}, 225--239.

\noindent \textsc{Steutel, F.W. {\rm and}  van Harn, K.} (2004).
\textit{Infinite divisibility of probability distributions on the
real line.} Monographs and Textbooks in Pure and Applied
Mathematics, 259. Marcel Dekker, New York.

 \noindent \textsc{Thorin, O.} (1977). On the
infinite divisibility of the lognormal distribution. \emph{cand.
Actuar. J.} \textbf{3}, 121--148.

\noindent \textsc{Thorin, O.} (1978) An extension of the notion of
a generalized $\Gamma $-convolution. \emph{Scand. Actuar. J.}
\textbf{3}, 141--149.



\end{footnotesize}

\bigskip

\begin{footnotesize}

\begin{tabbing} \= \textsc{Lancelot F. James}\\
\> \textsc{Department of Information and Systems Management}\\
\>\textsc{Hong Kong University of Science and Technology}\\
\> \textsc{Clear Water Bay, Kowloon, Hong Kong}\\
\> \textsc{E-mail:} \= \texttt{lancelot@ust.hk}
 \end{tabbing}

\medskip

 $ $ \hfill     \begin{tabular}{l}
    \textsc{Antonio Lijoi, Igor Pr\"unster}\\
        \textsc{Dipartimento di Economia Politica e Metodi Quantitativi}\\
        \textsc{Universit\`a degli Studi di Pavia}\\
        \textsc{Via San Felice 7, 27100 Pavia}\\
        \textsc{E-mail:} \texttt{lijoi@unipv.it, igor.pruenster@unipv.it}
        \end{tabular}

\end{footnotesize}
\end{document}